\pdfoutput=1
\documentclass[11pt]{article}
\usepackage[square,authoryear]{natbib}
\usepackage{marsden_article}
\usepackage{amscd}

\input xy
\xyoption{all}

\title{Hamilton-Pontryagin Integrators on Lie Groups Part I: \\
Introduction \& Structure-Preserving Properties}

\author{Nawaf Bou-Rabee\thanks{Applied
\& Computational Math., Caltech, Pasadena, CA 91125 ({\tt
nawaf@acm.caltech.edu}).}  \and
Jerrold E.~Marsden\thanks{Control \& Dynamical Systems, Caltech,
Pasadena, CA 91125 ({\tt marsden@cds.caltech.edu}). Research partially
supported by the National Science Foundation through NSF grant DMS-0204474.}}

\date{December, 2007}

\begin{document}

\maketitle

\begin{abstract} 
In this paper structure-preserving time-integrators for rigid body-type mechanical systems are derived from a discrete Hamilton-Pontryagin variational principle.  From this principle one can derive a novel class of variational partitioned Runge-Kutta methods on Lie groups.  Included among these integrators are generalizations of symplectic Euler and St\"{o}rmer-Verlet integrators from flat spaces to Lie groups.  Because of their variational design, these integrators preserve a discrete momentum map (in the presence of symmetry) and a symplectic form.  

In a companion paper, we perform a numerical analysis of these methods and report on numerical experiments on the rigid body and chaotic dynamics of an underwater vehicle.  The numerics reveal that these variational integrators possess structure-preserving properties that methods designed to preserve momentum (using the coadjoint action of the Lie group) and energy (for example, by projection) lack.  
\end{abstract}

\tableofcontents

\section{Introduction}

\paragraph{Overview.} This paper is concerned with efficient, structure-preserving time integrators for mechanical systems whose configuration space is a Lie group based on the Hamilton-Pontryagin (HP) variational principle \cite{Livens1919, LaWe2006, KWTKMSD2006, YoMa2006a, YoMa2006b}.  This HP Principle has many attractive theoretical properties; for instance, how it handles degenerate Lagrangian systems. The present paper paper shows that the HP viewpoint also provides a practical way to design discrete Lagrangians, which are the cornerstone of variational integration theory.   This overview explains the central idea of this paper in the context of vector spaces and shows how this idea extends to Lie groups.

The HP principle states that a mechanical system traverses a path that extremizes the following HP action integral: 
\begin{equation}
\mathfrak{S}_{HP} =  \underset{\text{Lagrangian}}{\underbrace{\int_a^b L(q, v) dt}} +  \underset{\text{kinematic constraint}}{\underbrace{\int_a^b \left\langle p, \dot{q} - v \right\rangle dt}} \text{.}
\end{equation}
The integrand of the HP action integral consists of two terms: the Lagrangian and a kinematic constraint paired with a Lagrange multiplier (the momentum).  The kinematic constraint relates the mechanical system's velocity to a curve on the tangent bundle. In this principle, the curves $q (t)$, $v (t)$, $p (t)$ are all varied independently. If $p $ is varied first, it collapses to the {\it usual Hamilton principle}. If, on the other hand, $v (t)$ is varied first it defines the (negative of the) Hamiltonian as the extrema of the terms involving $v$ and then the principle reduces to {\it Hamilton's phase space principle}.  This HP form of the action integral makes it amenable to discretization.   

In particular, one can implement an $s$-stage Runge-Kutta (RK) discretization of the kinematic constraint and enforce this discretization as a constraint in a discrete action sum.   The motivation is that the theory, order conditions, and implementation of such methods, are mature.  For this purpose let $[a,b]$ and $N$ be given, and define the fixed step size $h=(b-a)/N$ and $t_k = h k + a$, $k=0,...,N$.   Let $s$ be the number of stages in the RK method.  In analogy with the continuous system, the discrete HP action sum takes the following form:
\begin{align}
\mathfrak{S}_{HP}^d &=  \underset{\text{discrete Lagrangian}}{\underbrace{\sum_{k=0}^{N-1} \sum_{i=1}^s h b_i L(Q^i_k, V^i_k)}} \nonumber \\
&+ \underset{\text{discrete kinematic constraint}}{\underbrace{\sum_{k=0}^{N-1} \sum_{i=1}^s h \left\langle p_k^i, \frac{Q_k^i - q_k}{h} - \sum_{j=1}^s a_{ij} V_k^j \right\rangle + h \left\langle p_k, \frac{q_{k+1} - q_k}{h} - \sum_{j=1}^s b_j V_k^j\right\rangle }} \text{.}
\end{align}
It consists of two parts: a weighted sum of the Lagrangian using the weights from the Butcher tableau of the RK scheme, and pairings between  discrete internal and external stage Lagrange multipliers and the discretized kinematic constraint.    This strategy yields a Lagrangian analog of a well-known class of symplectic partitioned Runge-Kutta methods including the Lobatto IIIA-IIIB pair which generalize to higher-order accuracy \cite{Su1990, MaWe2001, HaLuWa2006}.

In the Lie group context, one can generalize this strategy using either constrained or generalized coordinates.  To use constrained coordinates one treats the system as a holonomically constrained mechanical system.  In this approach one assumes that $G$ can be written as the level set of some function $g: \mathbb{R}^n \to \mathbb{R}^k$, embeds $G$ in a larger linear space, and uses Lagrange multipliers to enforce the constraint.  This approach is discussed in \cite{BoOw2007b}.  The corresponding constrained action takes the following form:
\begin{align}
\mathfrak{S}_d^c = \sum_{k=0}^{N-1} \sum_{i=1}^{s} & h \left[ b_i L(Q_k^i, V_k^i)  + \left\langle p_k^i, \frac{Q_k^i - q_k}{h} -   \sum_{j=1}^s a_{ij} V_k^j \right\rangle  \right.  \nonumber \\
&+ \left. \left\langle p_{k+1}, \frac{q_{k+1} - q_k}{h} -   \sum_{j=1}^s b_{j} V_k^j \right\rangle +  b_i \left\langle \Lambda_k^i, g(Q_k^i) \right\rangle \right]  \text{.}
\end{align}

In the present paper a second approach based on generalized coordinates is presented.  First the paper introduces the following left-trivialized action:
\begin{equation}
\mathfrak{s}_{HP} =  \underset{\text{left-trivialized Lagrangian}}{\underbrace{\int_a^b \ell(g, \xi) dt}} +  \underset{\text{reconstruction equation}}{\underbrace{\int_a^b \left\langle \mu, g^{-1} \dot{g} - \xi \right\rangle dt}} \text{.}
\end{equation}
Then an equivalence is established between critical points of $\mathfrak{s}_{HP}$ and $\mathfrak{S}_{HP}$.  If the Lagrangian is left-invariant, it is shown that this principle unifies the system's Lie-Poisson and Euler-Poincar\'e descriptions  \cite{MaSc1993, CeMaPeRa2003}.  Since the reconstruction equation is a differential equation on a Lie group, one cannot directly discretize it by an RK method.  However, one can discretize it using an $s$-stage Runge-Kutta-Munthe-Kaas (RKMK) method \cite{Mu1995, MuZa1997, Mu1998, MuOw1999}.  The integral of the left-trivialized Lagrangian is approximated using a weighted sum given by the $b$-vector in the Butcher tableau of the RKMK scheme.   This approach is shown to yield a novel class of variational partitioned Runge-Kutta (VPRK) methods on Lie groups; including generalizations of  symplectic Euler and St\"{o}rmer-Verlet methods on flat spaces.

\section{Background and Setting}

In the next paragraphs we will give some background material for the reader's convenience as well as to put the paper into context.

\medskip

\noindent {\bf Variational Integrators.} Variational integration theory derives integrators for mechanical systems from discrete variational principles.   The theory includes discrete analogs of the Lagrangian, Noether's theorem, the Euler-Lagrange equations, and the Legendre transform.  Variational integrators can readily incorporate holonomic constraints (via Lagrange multipliers or the discrete null-space method; \cite{LeMaOr2007}) and non-conservative effects (via their virtual work) \cite{MaWe2001}, as well as discrete optimal control (see \cite{LeObMaOr2007} and references therein).    Altogether, this description of mechanics stands as a self-contained theory of mechanics akin to Hamiltonian, Lagrangian or Newtonian mechanics.

\begin{figure}[ht!]
\begin{center}
\includegraphics[scale=0.25,angle=0]{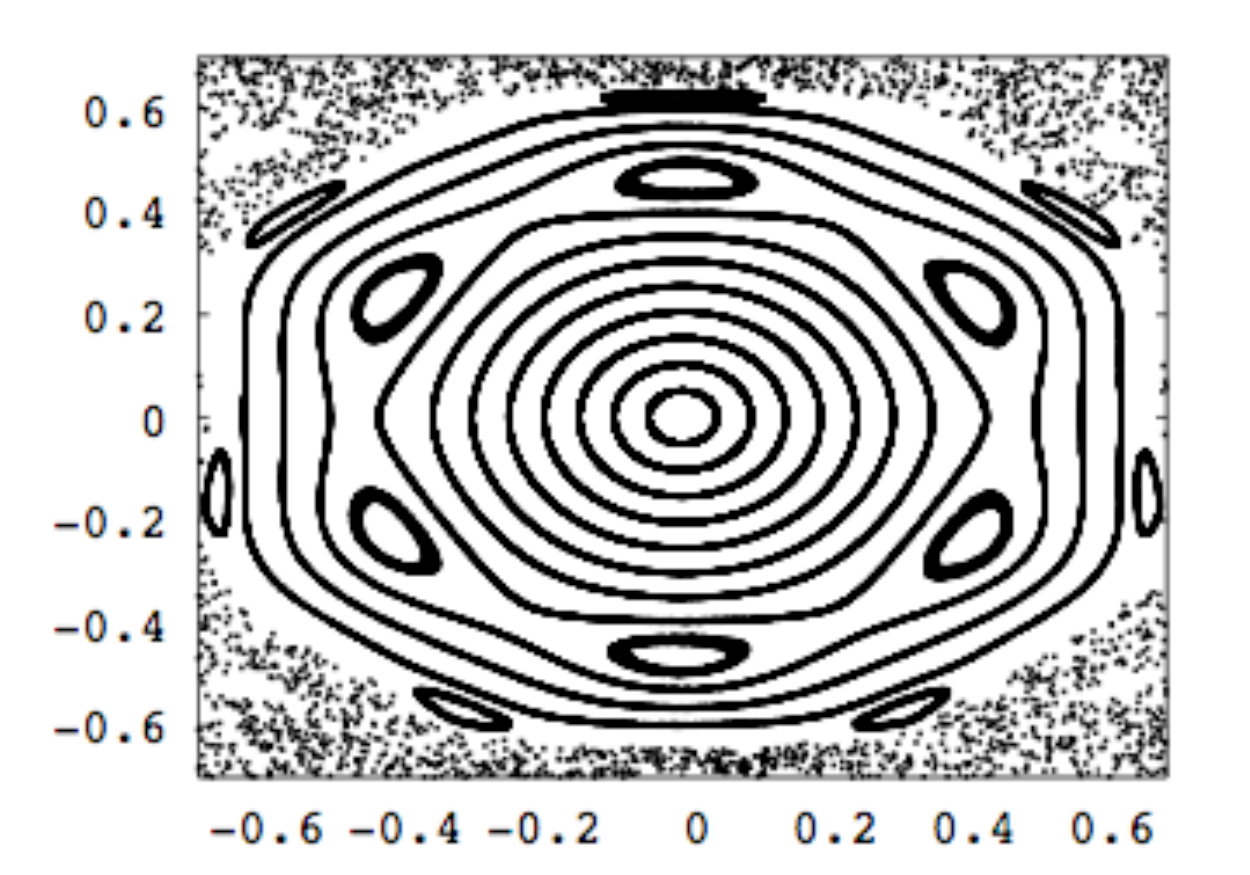}
\includegraphics[scale=0.25,angle=0]{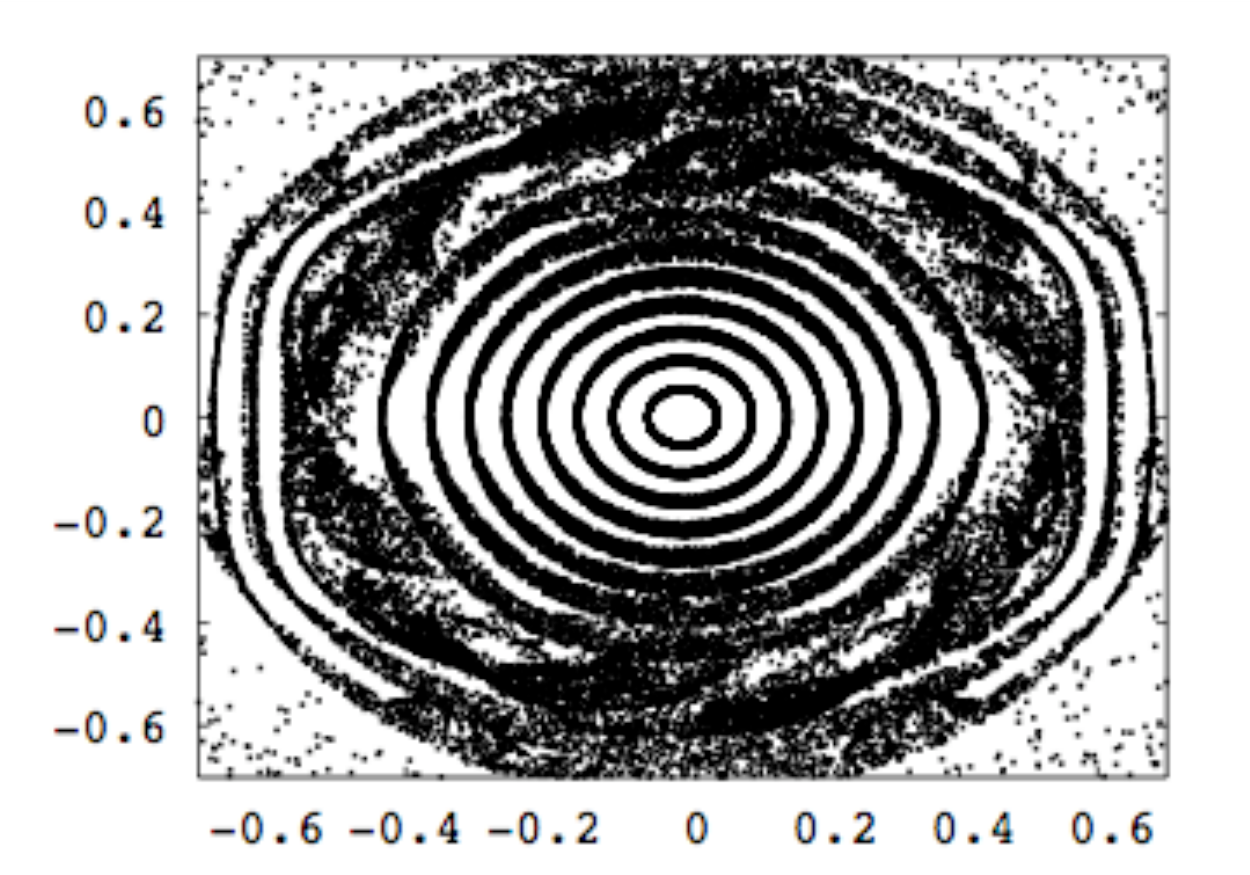} 
\includegraphics[scale=0.25,angle=0]{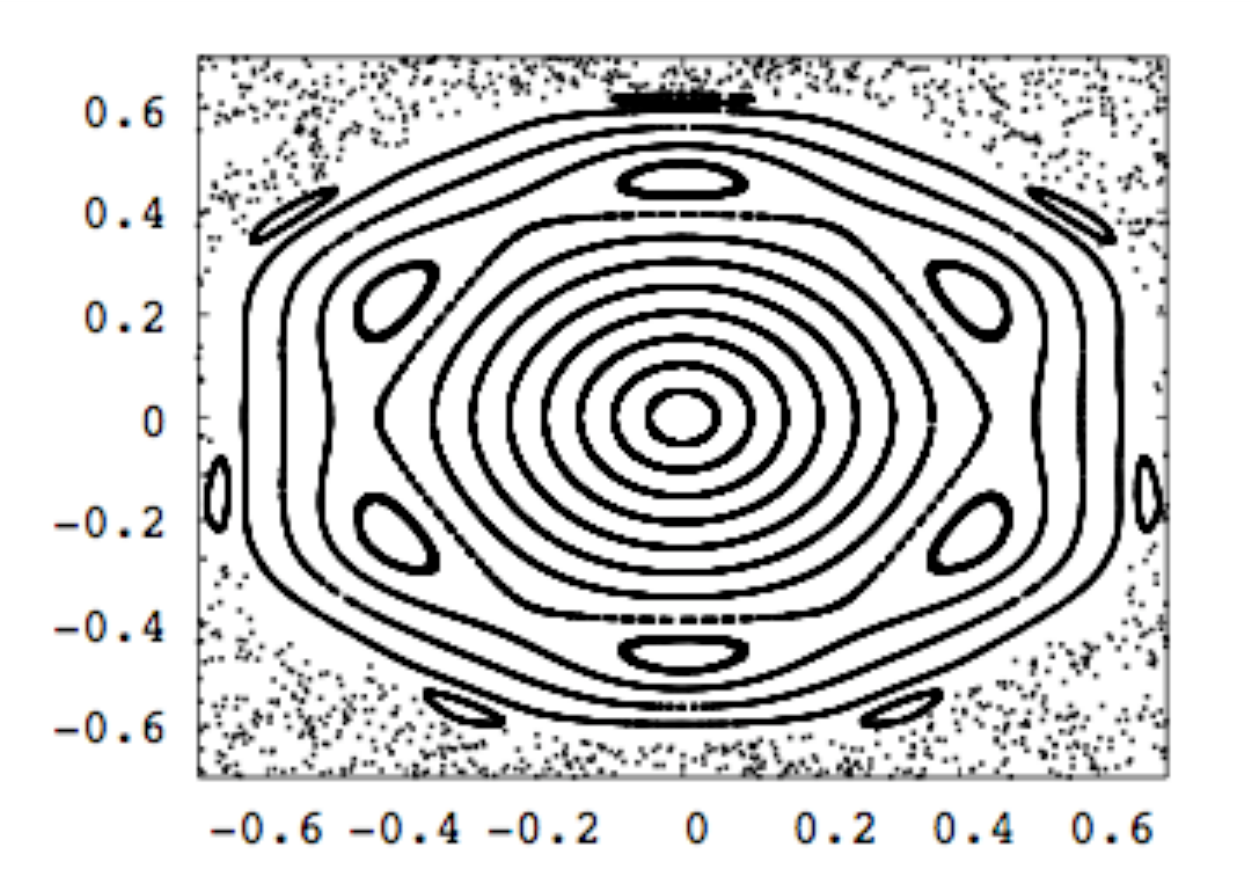}
\includegraphics[scale=0.25,angle=0]{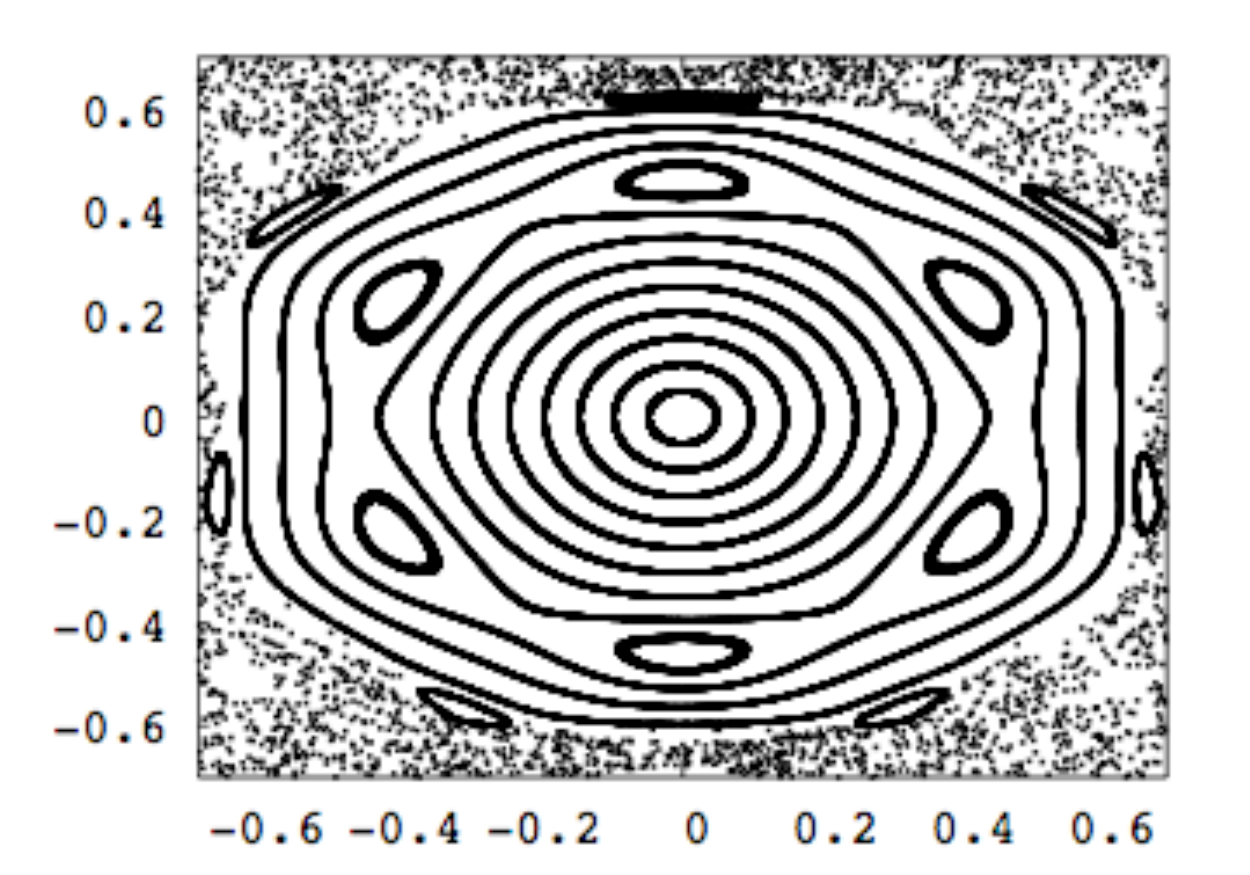} 
\hbox{\hspace{0.0in}   $\begin{array}{c} \text{(a) RK4} \\ $h=0.025$ \end{array}$ \hspace{0.6in}  $\begin{array}{c} \text{(b) RK4} \\ $h=0.05$ \end{array}$ \hspace{0.6in}  $\begin{array}{c} \text{(c) VE} \\ $h=0.025$ \end{array}$ \hspace{0.6in}  $\begin{array}{c} \text{(d) VE} \\ $h=0.05$ \end{array}$ } 
\caption{ \small {\bf Underwater Vehicle Dynamics.} 
This figure shows a computation of Poincar\'{e} sections using a second-order accurate variational Euler integrator (VE)  as compared to fourth-order accurate Runge-Kutta (RK4).  Both methods agree with the benchmark at the finer stepsize $h=0.025$.  However, at the coarser stepsize $h=0.05$, RK4 corrupts chaotic invariant sets while the lower-order accurate VE method preserves the structure of the benchmark.    } \label{fig:uv} \end{center}
\end{figure}

One of the distinguishing features of variational integrators is their ability to compute statistical properties of mechanical systems, such as in computing Poincar\'e sections, the instantaneous temperature of a system, etc.  For example, as a consequence of their variational design, variational integrators are symplectic.   A single-step integrator applied to a mechanical system is called {\em symplectic} if the discrete flow map it defines exactly preserves the canonical symplectic 2-form and is otherwise called standard.  Using backward error analysis one can show that symplectic integrators applied to Hamiltonian systems nearly preserve the energy of the continuous mechanical system for exponentially long periods of time and that the modified equations are also Hamiltonian \cite{HaLuWa2006}.  Standard integrators often introduce spurious dynamics in long-time simulations, e.g., artificially  corrupt chaotic invariant sets is well--illustrated in a computation from \cite{Bou-Rabee2007}, namely of a Poincar\'{e} section of an underwater vehicle in Fig.~\ref{fig:uv}  using a fourth-order accurate Runge-Kutta (RK4) method and a variational Euler (VE) method designed for rigid-body type systems.


In addition to correctly computing chaotic invariant sets and long-time excellent energy behavior, evidence is mounting that variational integrators correctly compute other statistical quantities in long-time simulations.  For example, in a simulation of a coupled spring-mass lattice, \cite{LeMaOrWe2004b} found that variational  integrators correctly compute the time-averaged instantaneous temperature (mean kinetic energy over all particles) over long-time intervals, whereas standard methods (even a higher-order accurate one) exhibit a artificial drift in this statistical quantity.   These structure-preserving properties of variational integrators motivated their extension to stochastic Hamiltonian systems.

\paragraph{Structure-Preserving Lie Group Integrators.} For a mechanical system on a Lie group that possesses the symmetry of that Lie group, in addition to the symplectic structure, the resulting flow preserves a momentum map associated with the Lie group symmetry.  In this context there are several different strategies available to derive structure-preserving Lie group integrators; some of these are discussed here.

One strategy involves the so-called Lie-Newmark method due to \cite{SiVu1988} and \cite{SiWo1991}.  These methods were motivated by the need to develop conserving algorithms that efficiently simulate the structural dynamics of rods and shells.   For example, the configuration space of a discrete, three-dimensional finite-strain rod model, would involve $N$ copies of $\mathbb{R}^3 \times \operatorname{SO(3)}$ where $N$ is the number of points in the discretization of the line of centroids of the rod.  For each point on the line of centroids, the orientation of the rod at that point is specified by an element of $\operatorname{SO}(3)$.  In such models the mathematical description of the rotational degrees of freedom at these points is equivalent to the EP description of a free rigid body with added nonconservative effects due to the elastic coupling between points.

It was not apparent  that the proposed Lie-Newmark methods had the necessary structure-preserving properties.  In fact, Simo \& Wong proposed another set of algorithms  which preserve momentum by using the coadjoint action on $\operatorname{SO}(3)$ to advance the flow.  Such integrators will be referred to as {\em coadjoint-preserving methods}.    Later, \cite{AuKrWa1993} showed that the midpoint rule member of the Lie-Newmark family with a Cayley reconstruction procedure was, in fact, a coadjoint-preserving method for $\operatorname{SO}(3)$.   They also numerically demonstrated the method's good performance crediting it to third-order accuracy in the discrete approximation to the Lie-Poisson structure.   In related work, \cite{McSc1995} construct reduced, coadjoint-orbit preserving integrators by reducing $G$-equivariant integrators on $T^*G$ obtained by embedding $G$ in a linear space using holonomic constraints.

Coadjoint and energy preserving methods of the Simo \& Wong type that further preserve the symplectic structure were developed for $\operatorname{SO}(3)$ by \cite{LeSi1994, LeSi1996}.  This was done by defining a one-parameter family of coadjoint and energy-preserving algorithms of the Simo \& Wong type in which the free parameter is a functional.  The function was specified so that the resulting map defined a transformation which preserves the continuous symplectic form.

Endowing coadjoint methods with energy-preserving properties was also the subject of \cite{EnFa2001}.  Specifically, they introduced integrators of the Runge-Kutta Munthe-Kaas type that preserved coadjoint orbits and energy using the coadjoint action on $\operatorname{SO}(3)$ and a numerical estimate of the gradient of the Hamiltonian.

Variational integration techniques have been used to derive structure-preserving integrators on Lie groups; see  \cite{MoVe1991, WeMa1997, MaPeSh1998, BoSu1999a, BoSu1999b}.   Moser and Veselov derived a variational integrator for the free rigid body by embedding $\operatorname{SO}(3)$ in the linear space of $3 \times 3$ matrices, $\mathbb{R}^9$, and using Lagrange multipliers to constrain the matrices to $\operatorname{SO}(3)$. This procedure was subsequently generalized to Lagrangian systems on more general configuration manifolds by the introduction of a discrete Hamilton's principle on the larger linear space with holonomic constraints to constrain to the configuration manifold in \cite{WeMa1997}.    They also considered the specific example of deriving a variational integrator for the free rigid body on the Lie group $S^3$ by embedding $S^3$ into $\mathbb{R}^4$  and using a holonomic constraint.  The constraint ensured that the configuration update remained on the space of unit quaternions (a Lie group) and  was enforced using a Lagrange multiplier.

Another approach is to use {\em reduction} to derive variational integrators on {\em reduced spaces}. \cite{MaPeSh1998}  developed a discrete analog  of EP reduction theory from which one could design {\em reduced} numerical algorithms.  They did this by constructing a discrete Lagrangian on $G \times G$ that inherited the $G$-symmetry of the continuous Lagrangian, and restricting it to the reduced space $(G \times G) / G \sim G$.  Using this discrete reduced Lagrangian and a discrete EP (DEP) principle, they derived DEP algorithms on the discrete reduced space. \label{acro:dep} They also considered using generalized coordinates to parametrize this discrete reduced space, specifically the exponential map from the  Lie algebra to the Lie group.  These techniques were applied to bodies with attitude-dependent potentials, discrete optimal control of rigid bodies, and to higher-order accuracy in  \cite{LeMcLe2005} and \cite{LeLeMc2007}.

\cite{BoSu1999a} considered a more general case where the symmetry group is a subgroup of the Lie group $G$ in the context of semidirect Euler-Poincar\'e theory (see \cite{HoMaRa1998}).  They did this by writing down the discrete Euler Lagrange equations for this system and left-trivializing them.  For the case when the symmetry group is $G$ itself, one recovers the DEP algorithm as pointed out in \cite{MaPeSh1998}.    In addition,  \cite{BoSu1999b} used this theory to determine and analyze an elegant, integrable discretization of the Lagrange top.

The perspective in this paper on Lie group variational integrators is different.   Recognizing that Euler's equations for a rigid body are in fact  decoupled from the dynamics on  the Lie group, and more generally, that the EP equation is decoupled from  the dynamics on the Lie group, the paper aims to develop discrete variational schemes that analogously consist of a reconstruction rule and discrete EP equations that can be solved independently of the reconstruction equation and on a lower dimensional linear space.  As mentioned in the overview the central idea is to discretize the reduced HP principle.

\paragraph{Organization of the Paper.} In \S3 continuous HP mechanics and its reduction is presented.  In particular, it is shown that the reduced and unreduced HP variational principle are equivalent to Hamilton's  and the EP variational principles.   Moreover, properties of the HP flow map are verified mainly to guide the discrete theory. In \S4 the reduced discrete analog of the HP theory is developed.   Properties of the discrete flow map are verified including discrete momentum map and symplectic form preservation. The theory is illustrated on several specific examples.   In \S5 the structure-preserving Lie group integrators relevant to this paper are presented.  In \S6 the free rigid body and underwater vehicle examples are presented, the structure-preserving methods from \S5 specialized to these examples, and results of numerical experiments are presented.  

\paragraph{Part II of this Paper.} The second installment of this paper will be devoted to the numerical analysis of HP methods along with  numerical experiments on a class of nonreversible mechanical systems on Lie groups as well as the chaotic dynamics of an underwater vehicle.   A specific outline of that paper is given in the conclusion section of the present paper.


\section{HP Mechanics} 

This section develops basic mechanics on Lie groups from the Hamilton-Pontryagin perspective.

\paragraph{The HP Principle.} Consider a mechanical system whose configuration space is a Lie group $G$. Let its tangent and cotangent bundles be denoted $TG$ and $T^*G$ respectively, and its Lie algebra and dual be given by $\mathfrak{g}$ and $\mathfrak{g}^*$ respectively.  In this paragraph the left-trivialization of the HP principle for a Lagrangian $L: TG \rightarrow \mathbb{R}$ will be derived.

The HP principle unifies the Hamiltonian and Lagrangian descriptions of a mechanical system, as shown in \cite{YoMa2006a, YoMa2006b}.   It states the following critical point condition on $TG \oplus T^*G$, 
\[
\delta \int _a^b \left[ L(g, v) + \left\langle p, \dot{g} - v \right\rangle \right] \,dt = 0,
\]  
where $(g(t),v(t),p(t)) \in TG \oplus T^*G$ are varied arbitrarily and independently with endpoint conditions $g (a) $ and $g (b)$ fixed. This builds in the Legendre transformation as well as the Euler--Lagrange equations into one principle.

\begin{definition}
Following standard conventions, the left action of $G$ on $TG$ or $T^*G$ is denoted by simple concatentation. The left-trivialized Lagrangian $\ell: G \times \mathfrak{g} \to \mathbb{R}$ is defined as:
\[
\ell(g, \xi) = L(g, g \xi)  \text{.}
\]
\end{definition}

The HP principle for mechanical systems on Lie groups is equivalent to the left trivialized HP principle: 
\[
\delta \int _a^b \left[ \ell (g, \xi) + \left\langle \mu, g  ^{-1} \dot{g} - \xi \right\rangle \right] \,dt = 0,
\]  
where there are no constraints on the variations; that is, the curves $\xi (t) \in \mathfrak{g}$, $\mu (t) \in \mathfrak{g}^\ast$ and $g (t) \in G $ can be varied arbitrarily. To see this, we proceed as follows.

Let $S(g,v,p)$ denote the HP action functional or integral,
\[ 
S(g,v,p) = \int _a^b \left[ L(g, v) + \left\langle p, \dot{g} - v \right\rangle \right] \,dt \text{.} 
\]   
Fixing the interval $[a,b]$, we regard $S$ as a map on {\it path space}: $S: \mathcal{C}(TG \oplus T^*G) \to \mathbb{R}$, where 
\[
\mathcal{C}(TG \oplus T^*G) = \{ (g,v,p): [a,b] \to TG \oplus T^*G ~|~ (g,v,p) \in C^{\infty}([a,b], TG \oplus T^*G) \}  \text{.}
\] 
Then a simple calculation shows that,
\begin{align*} 
S(g, v,p) & = \int_a^b \left[ L( g, g \xi) + \left\langle p,  g  g^{-1} (\dot{g} - v  ) \right\rangle \right] \, dt  \\
&= \int_a^b \left[ \ell(g, \xi) + \left\langle g p,  g^{-1} (\dot{g} - v  ) \right\rangle \right] \, dt \\
&=  \int_a^b \left[  \ell(g, \xi) + \left\langle \mu, g^{-1} \dot{g} - \xi \right\rangle \right] \, dt \\
&= s(g, \xi, \mu)
\end{align*} 
where $s$ is the reduced HP action functional, $\xi = g ^{-1} v \in \mathfrak{g}$,
and $\mu = g ^{-1} p \in \mathfrak{g}^*$.  From this equality one can derive the 
following key theorem.

\begin{theorem} Consider a Lagrangian system on a Lie group $G$ with Lagrangian $L: TG \to \mathbb{R}$.  Let $\ell:G \times \mathfrak{g} \to \mathbb{R}$ be its left-trivialization.   Then the following are equivalent
\begin{enumerate}
\item Hamilton's principle for $L$ on $G$ 
\[ 
\delta \int_a^b L(g, \dot{g}) dt =0 
\]
holds, for arbitrary variations $g(t)$ with  endpoint conditions $g (a) $ and $g (b)$ fixed;

\item the following variational principle holds on $\mathfrak{g}$,
\[ 
\delta \int_a^b \ell(g, \xi )dt = 0 
\] 
using variations of the form
\[ 
\delta \xi = \dot{\eta} + \operatorname{ad}_{\xi} \eta 
\]
where $\eta(a)=\eta(b)=0$ and $\xi = g^{-1} \dot{g}$; i.e., $\xi=T L_{g^{-1}} \dot{g}$;

\item the HP principle 
\[
\delta \int _a^b \left[ L(g, v) + \left\langle p, \dot{g} - v \right\rangle \right] \,dt = 0
\] 
holds, where $(g(t),v(t),p(t)) \in TG \oplus T^*G$, can be varied arbitrarily and independently  with endpoint conditions $g (a) $ and $g (b)$ fixed;

\item the left-trivialized HP principle 
\[
\delta \int _a^b \left[ \ell (g,\xi) + \left\langle \mu, g  ^{-1} \dot{g} - \xi \right\rangle \right] \,dt = 0 
\] 
holds, where $(g(t),\xi(t), \mu(t)) \in G \times \mathfrak{g} \times \mathfrak{g}^*$ can be varied arbitrarily and independently with endpoint conditions $g (a) $ and $g (b)$ fixed.
\end{enumerate}
\label{thm:hp}
\end{theorem}

\paragraph{Remark.}
If the Lagrangian is left-invariant, i.e., $L(g, v) = L(h g , h v)$ for all $h \in G$,  then the left-trivialized Lagrangian simplifies.  In particular, taking $h = g ^{-1}$, $\ell(\xi) = L(g, g \xi) = L(e, \xi )$, where $e$ is the identity element of the group.   In this case the left-trivialized HP principle unifies the Euler-Poincar\'e and Lie-Poisson descriptions on $\mathfrak{g}$ and $\mathfrak{g}^*$ respectively, consistent with the results of \cite{MaSc1993} and \cite{CeMaPeRa2003}.  

\paragraph{The HP Flow.} From the left-trivialized HP principle, the variations of $s$ with respect to $\xi$ and $\mu$ give 
\begin{align}
\textup{varying} \;  \mu \; \textup{gives} \quad  \xi = g^{-1} \dot{g}  \label{eq:rhpa} \text{,} ~~~ &\text{(reconstruction equation)} \text{,} \\
\textup{varying}\; \xi \; \textup{gives} \quad  \mu = \frac{\partial \ell}{\partial \xi}(g, \xi) \label{eq:rhpb} \text{,} ~~~ &\text{(Legendre transform)} \text{.}
\end{align}  
Also, setting the variation of $s$ with respect to $g$ equal to zero gives
 \begin{align}
& \int _a^b  \left[  \left\langle \frac{\partial \ell}{\partial g}, \delta g \right\rangle + \left\langle \mu, \delta (g^{-1} \dot{g} ) \right\rangle \right] \, dt \nonumber \\
& \qquad =  \int _a^b \left[  \left\langle g \frac{\partial \ell}{\partial g}, g^{-1} \delta g \right\rangle +  \left\langle \mu, - g^{-1} \delta g g^{-1} \dot{g} + g^{-1} \delta \dot{g} \right\rangle \right] \, dt = 0 \label{varyg}
\end{align}  
Observe that
 \[
\int _a^b \left[ \left\langle \mu, \delta (g^{-1} \dot{g} ) \right\rangle \right] \, dt = 
\int _a^b \left[  \left\langle \mu, - g^{-1} \delta g g^{-1} \dot{g} + g^{-1} \delta \dot{g} \right\rangle \right] \, dt 
\]   
Let $\eta =g^{-1} \delta g$. Using the product rule and \eqref{eq:rhpa}, we see that
\[
\frac{d}{dt} \eta = - \xi \eta + g^{-1} \frac{d}{dt} \delta g, \; \textup{which implies} \;
g^{-1} \frac{d}{dt} \delta g = \frac{d}{dt} \eta  + \xi \eta \text{.}
\]  
Substituting this relation into \eqref{varyg} gives 
\[
\int _a^b \left[ \left\langle g \frac{\partial \ell}{\partial g}, \eta \right\rangle +  \left\langle \mu, \frac{d}{dt}  \eta  + \operatorname{ad}_{\xi} \eta \right\rangle \right] \, dt = 0 \text{.}  
\]  
Integration by parts and using the boundary conditions on $g$ yields 
\[
\int _a^b \left[  \left\langle - \frac{d}{dt} \mu +  \operatorname{ad}_{\xi}^* \mu + g \frac{\partial \ell}{\partial g},   \eta  \right\rangle \right] \, dt =0 \text{.}
\]  
Since the variations are arbitrary, one arrives at
\begin{equation} 
\frac{d}{dt} \mu = \operatorname{ad}^*_{\xi} \mu + g \frac{\partial \ell}{\partial g} \label{eq:rhpc}  \text{.}
\end{equation}
In sum, the left-trivialized HP equations are given by:
\begin{equation} \label{eq:lthp}
\begin{cases}
\begin{array}{rcl}
\frac{d}{dt} g &=& g \xi  \text{,} \\
\frac{d}{dt} \mu &=& \operatorname{ad}^*_{\xi} \mu + g \frac{\partial \ell}{\partial g}  \text{,} \\
\mu &=& \frac{\partial \ell}{\partial \xi}(g, \xi) \text{.} 
\end{array}
\end{cases}
\end{equation}  
Assuming that the Legendre transform is invertible, \eqref{eq:lthp} describes an IVP on the left-trivialized space $G \times \mathfrak{g} \times \mathfrak{g}^*$.

\begin{definition}
Let $\mathcal{I}_{HP}$ denote the {\bfi admissible space} and defined as, 
\begin{equation}
\mathcal{I}_{HP} :=  \left\{ (g, v, p) \in TG \oplus T^*G \; \left| \;  p = \frac{\partial L}{\partial v}(g, v)  \right\} \right. \text{.} 
\end{equation}
Let $\mathcal{I}_{hp}$ denote its left-trivialization and defined as the subset of $G \times \mathfrak{g} \times \mathfrak{g}^*$ that satisfies \textup{\eqref{eq:rhpb}}, i.e., 
\begin{equation}
\mathcal{I}_{hp} :=  \left\{ (g, \xi, \mu) \in G \times \mathfrak{g} \times \mathfrak{g}^* \; \left| \;   \mu = \frac{\partial \ell}{\partial \xi}(g, \xi)  \right\} \right. \text{.} \label{eq:Ihp}
\end{equation}
The {\it natural projection} is denoted by $\pi_{HP}: TG \oplus T^*G \to T^*G$ and defined as, 
\[
\pi_{HP}(g, v, p) := (g, p),~~~\pi_{HP}^{-1}(g,p) = (g, v, p),~~~(g, v) = \mathbb{F}L^{-1}(g,p) 
\]   
where $\mathbb{F}L$ is the Legendre transform.
\end{definition}

Given a time-interval $[a,b]$ and an initial $(g(a), \xi(a), \mu(a)) \in \mathcal{I}_{hp}$, one can solve for $(g(b), \xi(b), \mu(b)) \in \mathcal{I}_{hp}$ by eliminating $\xi$ using the left-trivialized Legendre transform \eqref{eq:rhpb} and solving the ODEs  \eqref{eq:rhpa} and \eqref{eq:rhpc} for $g$ and $\mu$.  Let this map on $\mathcal{I}_{hp}$ be called the {\it left-trivialized HP flow map}, $F_{hp}: \mathcal{I}_{hp}  \to \mathcal{I}_{hp}$.

The flow map $F_{hp}$ is equivalent to the HP flow on $\mathcal{I}_{HP}$  through left trivialization which defines a diffeomorphism between $TG \oplus T^*G$ and $G \times \mathfrak{g} \times \mathfrak{g}^*$, and hence, between $\mathcal{I}_{HP}$ and $\mathcal{I}_{hp}$.   Through $\pi_{HP}$ the HP flow is identical to the Hamiltonian flow for the Hamiltonian of this mechanical system on $T^*G$ obtained via the Legendre transformation.  Although $\pi_{HP}$ is not a diffeomorphism from $TG \oplus T^*G$ to $T^*G$, it is a diffeomorphism when its domain is restricted to $\mathcal{I}_{HP}$.  Thus, the left-trivialized HP, HP and Hamiltonian flows of this mechanical system are all equivalent. This observation makes the subsequent proof of symplecticity seem superfluous, since this structure obviously follows from the standard  theory of Hamiltonian systems with symmetry.   However, this verification is still important since it serves as a model for the less obvious discrete theory.

It will be helpful to define $\pi_{\mathcal{I}_{HP}} = \pi_{HP} |_{\mathcal{I}_{HP}}$.  The manifold $TG \oplus T^*G$ is a presymplectic manifold with the {\em HP presymplectic form}, $\Omega_{HP} = \pi_{HP}^* \Omega$, and the manifold $\mathcal{I}_{HP}$ is a symplectic manifold with the {\em HP symplectic form}, $\Omega_{\mathcal{I}_{HP}}=  \pi_{\mathcal{I}_{HP}}^* \Omega$.  
Similarly, the manifold $G \times \mathfrak{g} \times \mathfrak{g}^*$ is a presymplectic manifold with the presymplectic form $\omega_{HP}$ that is obtained by pulling-back the HP presymplectic form by the left trivialization of $TG \oplus T^*G$, $\phi: G \times \mathfrak{g} \times \mathfrak{g}^* \to TG \oplus T^*G$, i.e., $\omega_{HP} = \phi^* \Omega_{HP}$.   However, if the left-trivialization is restricted to $\mathcal{I}_{hp}$,  $\phi_{\mathcal{I}_{hp}} = \phi |_{\mathcal{I}_{hp}}$, then $\mathcal{I}_{hp}$ is a symplectic manifold with the symplectic form given by $\omega_{\mathcal{I}_{hp}} = \phi_{\mathcal{I}_{hp}}^* \Omega_{\mathcal{I}_{HP}}$.

\paragraph{Symplecticity.}
The symplectic structure of left-trivialized HP flows is obvious from the standard theory of Hamiltonian systems with symmetry, but reviewing the proof will help since it parallels the discrete case.

Consider the restriction of the left-trivialized HP action integral to solutions of \eqref{eq:lthp}: $\hat{s}$.  Since the space of solutions of \eqref{eq:lthp} can be identified with $\mathcal{I}_{hp}$, $\hat{s}: \mathcal{I}_{hp} \to \mathbb{R}$.   The differential of $\hat{s}$ can be written as, 
\begin{align*}
\mathbf{d} \hat{s} \cdot (\delta g(a), \delta \xi(a), \delta \mu(a)) =& 
\int_a^b  \left[ \left( g^{-1} \dot{g} - \xi \right) \cdot \delta \mu +
\left( \mu - \frac{\partial \ell}{\partial \xi} \right) \cdot \delta \xi \right] \, dt \\
&+ \int_a^b\left[ \left( - \frac{d}{dt} \mu + \operatorname{ad}_{\xi}^* \mu + g \frac{\partial \ell}{\partial g} \right) \cdot g^{-1} \delta g \right] \, dt  
+  \left. \left\langle \mu, g^{-1} \delta g \right\rangle \right|_a^b  \\
=&  \left. \left\langle \mu, g^{-1} \delta g \right\rangle \right|_a^b = 
  ( (F_{hp})^* \theta_{\mathcal{I}_{hp}} - \theta_{\mathcal{I}_{hp}} ) \cdot (\delta g(a), \delta \xi(a), \delta \mu(a))
\end{align*}  
where we have introduced the left-trivialized HP one-form, $\theta_{\mathcal{I}_{hp}} = \phi_{\mathcal{I}_{hp}}^* \Theta_{\mathcal{I}_{HP}}$.  Since $\mathbf{d}^2 \hat{s} = 0$, observe that
\[ 
\mathbf{d}^2 \hat{s} = (F_{hp})^* \omega_{\mathcal{I}_{hp}} - \omega_{\mathcal{I}_{hp}}  = 0 \text{.}  
\]
And hence, as a map on $\mathcal{I}_{hp}$, $F_{hp}$ is symplectic.

\begin{theorem} Left-trivialized HP flows preserve the symplectic two-form $\omega_{\mathcal{I}_{hp}}$.  \end{theorem}


\section{Lie Group VPRK Integrators}

The purpose of this section is to use the general HP methodology to derive a variety of integrators of variational partitioned Runge-Kutta (VPRK) type for Lie groups. After introducing the map $\tau$ which is typified by the exponential map, and its properties, we use an $s$-stage  Runge-Kutta-Munthe-Kaas (RKMK) approximation to the reconstruction equation, which leads naturally to the introduction of VPRK Integrators on Lie groups. This includes the St\"{o}rmer-Verlet method for Lie groups, variational Euler methods on Lie groups, and Euler-Poincar\'e integrators.

\paragraph{Canonical Coordinates of the First Kind.} To setup the discrete HP principle, we introduce a map $\tau: \mathfrak{g} \to G$.  Let $e\in G$ be the identity element of the group. The map $\tau$ is assumed to be a local diffeomorphism mapping a neighborhood of  zero on $\mathfrak{g}$ to one of $e$ on $G$ with $\tau(0) = e$, assumed to be analytic in this neighborhood, and assumed to satisfy $\tau(\xi) \cdot \tau(-\xi) = e$.   Thereby $\tau$ provides a local chart on the Lie group.  By left translation this map can be used to construct an atlas on $G$.   An example of a $\tau$ is the exponential map on $G$, but there are other interesting examples as well, as we shall see shortly.

\begin{definition}
 The local coordinates associated with the map $\tau$ are called {\bfi canonical coordinates of the first kind} or just {\bfi canonical coordinates}.  
 \end{definition}

For an exposition of canonical coordinates of the first and second kind, and their applications  the reader is referred to \cite{IsMuNoZa2000}.   In what follows we will prove some properties of these coordinates that will be needed shortly.

\paragraph{Derivative of $\tau$ and its inverse.} To derive the integrator that comes from a discrete left-trivialized HP principle, we will need to differentiate $\tau^{-1}$.    The {\em right trivialized tangent} of $\tau$ and its inverse will play an important role in writing this derivative in an efficient way.  The following is taken from Definition 2.19 in \cite{IsMuNoZa2000}.

\begin{definition}  
Given a local diffeomorphism $\tau: \mathfrak{g} \to G$, we define its {\bfi right trivialized tangent} to be the function $d \tau: \mathfrak{g} \times \mathfrak{g} \to \mathfrak{g}$ which satisifies, 
\[
\operatorname{D} \tau( \xi )  \cdot \delta = TR_{\tau( \xi )} d \tau_{\xi} (\delta )  \text{.}
\]  
The function $d \tau$ is linear in its second argument. 
\label{def:dtau}
\end{definition}

\begin{figure}[ht]
\begin{center}
\includegraphics[scale=0.35,angle=0]{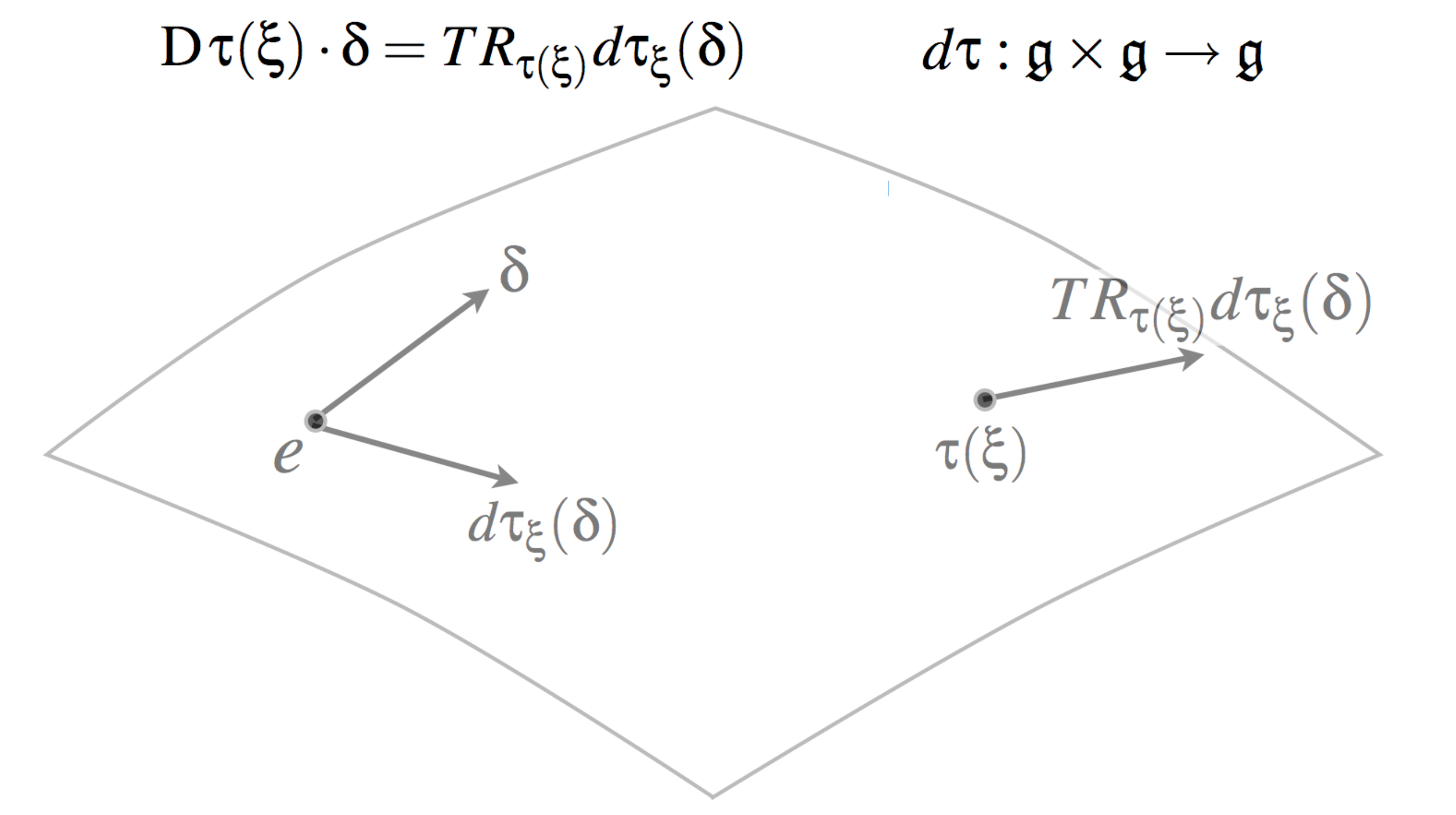}
\caption{\footnotesize  {\bf Derivative of $\boldsymbol{\tau}$.}  Definition \eqref{def:dtau} splits the differential of $\tau$ into a map on the Lie algebra (the right trivialized tangent of $\tau$) 
and right multiplication to the tangent space at $\tau(\xi)$.  }
\label{fig:dtau}
\end{center}
\end{figure}

Figure \ref{fig:dtau} illustrates the geometry behind this definition.

From this definition the following lemma is deduced.

\begin{lemma}  
The following identity holds, 
\[
d \tau_{\xi} (\delta )  = \operatorname{Ad}_{\tau(\xi)} d \tau_{- \xi} (\delta )  \text{.}
\]  
\label{dtaupropertya}
\end{lemma}

\begin{proof}
Differentiation of $\tau(\xi) \cdot \tau(-\xi) = e$ gives 
\[
 \operatorname{D}\tau( -\xi) \cdot \delta  = - TL_{\tau( - \xi)} TR_{\tau( -\xi)}  \left( \operatorname{D} \tau( \xi) \cdot \delta \right)  \text{.}
\]  
While the chain rule yields
\[
\operatorname{D}\tau (- \xi) \cdot \delta = - T R_{\tau( - \xi )} d \tau_{- \xi} (\delta )  \text{.}
\]  
Combining these two identities and using the definition above, 
\[
- T R_{\tau (-\xi)}  d \tau_{- \xi} (\delta ) =  - TL_{\tau( - \xi)} T R_{ \tau(-\xi) } TR_{\tau(\xi)}   d \tau_{\xi} (\delta )    \text{.}
\] 
Simplifying this expression gives, 
\[
T L_{\tau( \xi)} d \tau_{- \xi} (\delta )  = T R_{\tau (\xi)} d \tau_{\xi} (\delta )  \text{,}
\]  
which proves the identity.  
\end{proof}

We will also need a simple expression for the differential of $\tau^{-1}$.

\begin{definition} 
The {\bfi inverse right trivialized tangent} of $\tau$ is the function $d \tau^{-1}: \mathfrak{g} \times \mathfrak{g}  \to \mathfrak{g}$ which satisifies for $g = \tau(\xi)$, 
\[
\operatorname{D}\tau^{-1}( g )  \cdot \delta = d \tau^{-1}_{\xi} ( TR_{\tau( - \xi )} \delta ),~~~d \tau^{-1}_{\xi} ( d \tau_{\xi} (\delta)) = \delta  \text{.}
\]  
The function $d \tau^{-1}$ is always linear in its second argument. 
\label{def:dtauinverse}
\end{definition}

\begin{figure}[ht]
\begin{center}
\includegraphics[scale=0.35,angle=0]{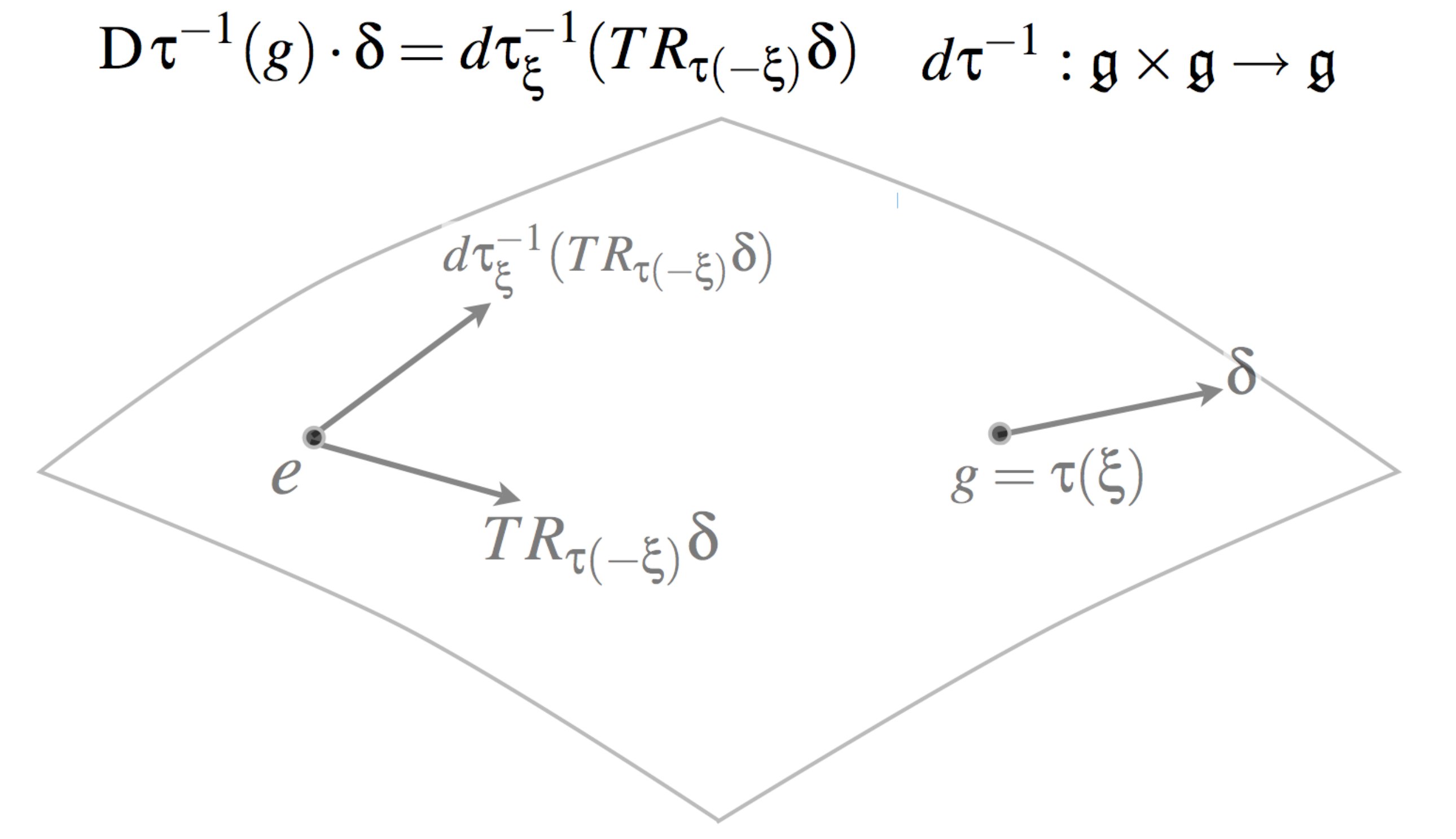}
\caption{\footnotesize  {\bf Derivative of $\boldsymbol{\tau}^{-1}$.}  Definition \ref{def:dtauinverse} splits the differential of $\tau^{-1}$ into right multiplication to the Lie algebra 
and a map on the Lie algebra (the right trivialized tangent of $\tau^{-1}$).  }
\label{fig:dtauinverse}
\end{center}
\end{figure}

Figure \ref{fig:dtauinverse} illustrates the geometry behind this definition.

The following lemma follows from this definition and Lemma~\ref{dtaupropertya} above.

\begin{lemma}  The following identity holds, 
\[
d \tau_{\xi}^{-1} (\delta )  =  d \tau_{- \xi}^{-1} (\operatorname{Ad}_{\tau(-\xi)} \delta )  \text{.}
\] 
\label{propertyb}
\end{lemma}

\begin{proof}
This follows directly from Lemma~\ref{dtaupropertya}.  Let $\delta \to d \tau_{\xi}^{-1} (\delta)$ in that identity to obtain  
\[
\delta  = \operatorname{Ad}_{\tau(\xi)} d \tau_{- \xi} (d \tau_{\xi}^{-1} (\delta) )  \text{.}
\]  
And now solve this equation for $d \tau_{\xi}^{-1}(\delta) $, 
\[
d \tau_{\xi}^{-1} (\delta)  = d \tau_{- \xi}^{-1} \left(\operatorname{Ad}_{\tau(-\xi)}  \delta \right) \text{.}
\] 
\end{proof}

\paragraph{RKMK Discretization of Reconstruction Equation.}

Let $[a,b]$ and $N$ be given, let $h = (b-a)/N$ be a fixed integration time step and $t _k = h k $.   A good candidate for discretizing the reconstruction equation is given by a generalization of $s$-stage Runge-Kutta methods to differential equations on Lie groups, namely Runge-Kutta-Munthe-Kaas (RKMK) methods introduced in the following series of papers: \cite{Mu1995, MuZa1997, Mu1998, MuOw1999}.  The idea behind those papers is to use canonical coordinates on the Lie group to transform the differential equation on $TG$, e.g., given by, 
\begin{equation}   \label{eq:liegroupdeq}
\dot{g} = g f(t,g), ~~~ g(0)=g_0,~~~ g(t)\in G,~~~f(t,g(t)) \in \mathfrak{g} \text{.}
\end{equation}
to a differential equation on $\mathfrak{g}$.  Specifically,  substitute the following parametrization $g(t) = g_0 \tau(\Theta(t)) $ into \eqref{eq:liegroupdeq} to obtain,
\[
\dot{g} = T L_{g_0} T R_{\tau(\Theta)} d \tau_{\Theta} \dot{\Theta} = T L_{g_0} TL_{\tau(\Theta)} f(t, g) \text{.}
\]
Using Lemma~\ref{dtaupropertya} this equation can be rewritten as,
\[
TL_{\tau(-\Theta)} T R_{\tau(\Theta)} d \tau_{\Theta} \dot{\Theta} = \operatorname{Ad}_{\tau(-\Theta)}  d \tau_{\Theta} \dot{\Theta} = d \tau_{-\Theta} \dot{\Theta}
= f(t, g) \text{.}
\]
Solving for $\dot{\Theta}$ gives
\begin{equation}
\dot{\Theta} = d \tau^{-1}_{-\Theta} f(t, g), ~~~\Theta(0) = 0,~~~\Theta(t) \in \mathfrak{g} \label{eq:liealgebradeq} \text{.}
\end{equation}
As described in the following definition, the RKMK method is obtained by applying an $s$-stage  RK method to \eqref{eq:liealgebradeq}.

\begin{definition}  
Consider the first-order differential equation $\dot{g} = f(t,g)$ for the curve $(g(t), f(t,g(t))) \in TG$.  Given coefficients $b_i, a_{ij} \in \mathbb{R}$ ($i, j=1,\cdots,s$) and set $c_i = \sum_{j=1}^s a_{ij}$.  An {\bfi $s$-stage  Runge-Kutta-Munthe-Kaas} (RKMK) approximation is given by 
\begin{align}
G_k^i &= g_k \tau(h \Theta_k^i ) \text{,} \label{eq:rkmkinternalA} \\
\Theta_k^i &=  h \sum_{j=1}^s  a_{ij} d\tau^{-1}_{-h \Theta_k^j} f(t_k+c_j h, G_k^j)  \text{,} ~~~ i=1, \cdots, s  \text{,} \label{eq:rkmkinternalB} \\
g_{k+1} &= g_k \tau\left(h \sum_{j=1}^s b_j d\tau^{-1}_{-h \Theta_k^j} f(t_k+c_j h, G_k^j) \right) \text{.} \label{eq:rkmkexternal}
\end{align}
If $a_{ij} = 0$ for $i \le j$ the RKMK method is called explicit, and implicit otherwise.   The vectors $g_k$ and $G_k^i$ are called external and internal stage configurations, respectively.  \label{acro:rkmk}
\end{definition}

It follows that for given $\tau$ an $s$-stage  RKMK method is determined by its $a$-matrix and $b$-vector which are typically displayed using the so-called Butcher tableau: 
\begin{center}
\begin{tabular}{c|ccc}
    $c_1$ &  $a_{11}$ & $\cdots$ & $a_{1s}$ \\
$\vdots$ & $\vdots$   &             & $\vdots$ \\
  $c_s$  & $a_{s1}$   &  $\cdots$  & $a_{ss}$ \\
\hline  
          & $b_1$       & $\cdots$   & $b_s$ 
\end{tabular}
\end{center}
Suppose that $\xi(t)$, $t\in[a,b]$, is given.  From this definition it is clear that an $s$-stage  RKMK method applied to $\dot{g} = g \xi$ can be written as:
\begin{equation} \label{eq:rkmk}
\begin{cases}
\begin{array}{rclc}
\tau^{-1}( g_k^{-1} G_k^i)/h &=& \sum_{j=1}^s a_{ij} d \tau^{-1}_{-h \Theta_k^j} \Xi_k^j  = \Theta_k^i \text{,} & i=1,\cdots,s, \\
\tau^{-1}(g_k^{-1} g_{k+1})/h &=&  \sum_{j=1}^s b_j d \tau^{-1}_{-h \Theta_k^j} \Xi_k^j  \text{.} 
\end{array}
\end{cases}
\end{equation}
where $\Xi_k^i = \xi(t_k + c_i h)$.    In practice one often truncates the series expansion of $ d \tau^{-1}_{-h \Theta_k^j}$.  The following theorem guides how to do this without degrading the order of accuracy \cite{HaLuWa2006}.

\begin{theorem} \label{thm:dtautruncation}
Given a qth order approximant to the exact exponential: $\tau: \mathfrak{g} \to G$.  If the underlying RK method is of order $p$ and the truncation index of $ d \tau^{-1}_{-\Theta}$ satisfies $q \ge p-2$ then the RKMK method is of order $p$.
\end{theorem}

\paragraph{VPRK Integrators on Lie Groups.} The discrete HP principle states that the discrete path the discrete system takes is one that extremizes a reduced action sum that will be introduced shortly.   To discretize the action integral, \eqref{eq:rkmk} is treated as a constraint in the discrete HP action, and the integral of the left-trivialized Lagrangian is approximated by the following quadrature:
\begin{equation}
\int_{t_k}^{t_k+h} \ell( g, \xi ) dt \approx \sum_{i=1}^s h b_i \ell(G_k^i, \Xi_k^i )  \text{.}
\end{equation}
The truncation index of $d \tau^{-1}_{-\Theta}$ in \eqref{eq:rkmk} is chosen to be $q=0$.  By theorem~\ref{thm:dtautruncation} one can obtain second-order accurate methods from this principle.

\begin{definition}
Given an $s$-stage  RKMK method with $b_j \ne 0$ for $j=1,...,s$, define the {\bfi discrete VPRK path space},  
\begin{align*}
\mathcal{C}_d(g_1,g_2) = \{ & (g, \mu, \{ \Theta^i, \Xi^i, \mu^i \}_{i=1}^s)_d : \{ t_k \}_{k=0}^N \to (G \times \mathfrak{g}^*) \times (\mathfrak{g} \times \mathfrak{g} \times \mathfrak{g}^* )^s ~|~ \\
& g(t_0) = g_1,~~g(t_N) = g_2 \}  \text{.}
\end{align*}
and the {\bfi action sum}  $s_d : \mathcal{C}_d(g_1,g_2) \to \mathbb{R}$ as  
\begin{align}
s_d =  \sum_{k=0}^{N-1} \sum_{i=1}^s h & \left[ b_i \ell(G_k^i, \Xi_k^i)  
  +  \left\langle \mu_k^i, \tau^{-1}(g_k^{-1} G_k^i)/h - \sum_{j=1}^s a_{ij} \Xi_k^j \right\rangle \right.  \nonumber \\
& + \left. \left\langle \mu_{k+1}, \tau^{-1}(g_k^{-1} g_{k+1})/h - \sum_{j=1}^s b_j \Xi_k^j \right\rangle \right] \text{.} \label{eq:hpsumfuture}
 \end{align} 
\end{definition}

Observe that $s_d$ is an approximation of the reduced HP action integral by numerical quadrature.  The definition of $\tau$ as a map from $\mathfrak{g}$ to $G$ ensures that the pairings  in the above sum arae well defined.   The discrete left-trivialized HP principle states that, 
\[
 \delta s_d = 0
 \] 
for arbitrary and independent variations of  the external stage vectors $( g_k, \mu_k ) \in G \times \mathfrak{g}^* $ and the internal stage vectors $( \Theta_k^i, \Xi_k^i, \Psi_k^i ) \in \mathfrak{g} \times \mathfrak{g} \times \mathfrak{g}^*$ for $i=1, \cdots, s$ and $k=0, \cdots, N$ subject to  fixed endpoint conditions on $\{ g_k \}_{k=0}^N$.

\begin{theorem}
Let $\ell: G \times \mathfrak{g} \to \mathbb{R}$ be a smooth, left-trivialized Lagrangian.  A discrete curve $c_d \in \mathcal{C}_d(g_1,g_2)$ satisfies the following VPRK scheme:
\begin{align} 
 \tau^{-1}(g_k^{-1} G_k^i)/h & =   \sum_{j=1}^s a_{ij} \Xi_k^j = \Theta_k^i, \label{eq:lievprk1} \\
 \tau^{-1}(g_k^{-1} g_{k+1})/h & =    \sum_{j=1}^s b_j \Xi_k^j = \xi_{k+1},   \\
 (d \tau^{-1}_{h \xi_{k+1}})^* M^i_{k} &=    (d \tau^{-1}_{-h \xi_k})^* \mu_k   \nonumber \\
 & \hspace{-0.45in} + h \sum_{j=1}^s \left( b_j (d \tau^{-1}_{h \Theta_k^j})^* - \frac{b_j a_{ji}}{b_i} (d \tau^{-1}_{h \xi_{k+1}})^* \right) (d \tau_{-h \Theta^j_k})^* G^j_k \frac{\partial \ell}{\partial g}(G^j_k, \Xi^j_k ),   \\
 (d \tau^{-1}_{h \xi_{k+1}})^* \mu_{k+1} &=    (d \tau^{-1}_{-h \xi_k})^* \mu_k  + h \sum_{j=1}^s b_j (d \tau^{-1}_{h \Theta^j_k})^* (d \tau_{-h \Theta^j_k})^* G^j_k \frac{\partial \ell}{\partial g}(G^j_k, \Xi^j_k ),    \\
M^i_{k} &=    \frac{\partial \ell}{\partial \xi}(G^i_k, \Xi^i_k )  \text{.}  \label{eq:lievprk}
\end{align}
for $i=1,\cdots,s$ and $k=0,\cdots,N-1$, if and only if it is a critical point of the function $s_d:  \mathcal{C}_d(g_1, g_2) \to \mathbb{R}$, that is, $\mathbf{d} s_d (c_d) = 0$.   Moreover, the discrete flow map defined by the above scheme, $F_h:  \mathcal{I}_{hp}  \to \mathcal{I}_{hp}$, preserves the symplectic form $\omega_{\mathcal{I}_{hp}}$.  
\end{theorem}

\begin{proof}
Set $\eta_k = g_k^{-1} \delta g_k$ and $H_k^i = {G_k^{i}}^{-1} \delta G_k^i$.    The differential of $s_d(c_d)$  in the direction $z = (\{ \delta g_k, \delta \mu_k \}, \{ \delta G_k^i, \delta \Xi_k^i, \delta \mu_k^i \}_{i=1}^s)$ is given by:
\begin{align*}
\mathbf{d} s_d  \cdot z & = \sum_{k=0}^{N-1} \sum_{i=1}^s h b_i G^i_k \frac{\partial \ell}{\partial g}(G^i_k, \Xi^i_k ) \cdot H^i_k 
+ h b_i \frac{\partial \ell}{\partial \xi}( G^i_k, \Xi^i_k ) \cdot \delta \Xi^i_k  \\
& \quad + h \left\langle \delta \mu^i_k, \tau^{-1} (g_k^{-1} G_k^i)/h - \sum_{j=1}^s a_{ij} \Xi^j_k \right\rangle  \\
& \quad + h \left\langle \delta \mu_{k+1}, \tau^{-1} (g_k^{-1} g_{k+1})/h - \sum_{j=1}^s b_j \Xi^j_k \right\rangle \\
& \quad + h \left\langle \mu^i_k,  - d \tau^{-1}_{h \Theta^i_k} \eta_k  / h + d \tau^{-1}_{-h \Theta^i_k} H^i_k / h - \sum_{j=1}^s a_{ij} \delta \Xi^j_k \right\rangle  \\
& \quad + h \left\langle \mu_{k+1},  - d \tau^{-1}_{h \xi_{k+1}} \eta_k  / h + d \tau^{-1}_{-h \xi_{k+1}} \eta_{k+1} / h - \sum_{j=1}^s b_j \delta \Xi^j_k \right\rangle 
\end{align*}
Collecting terms with the same variations and summation by parts using the boundary conditions $\delta g_0 = \delta g_N = 0$ gives,
\begin{align*}
\mathbf{d} s_d  \cdot z & = \sum_{k=1}^{N-1}  \sum_{i=1}^s h \left\langle \delta \mu^i_k, \tau^{-1} (g_k^{-1} G_k^i)/h - \sum_{j=1}^s a_{ij} \Xi^j_k \right\rangle  \\
& \quad + h \left\langle \delta \mu_{k+1}, \tau^{-1} (g_k^{-1} g_{k+1})/h - \sum_{j=1}^s b_j \Xi^j_k \right\rangle  \\
& \quad + h \left\langle b_i \frac{\partial \ell}{\partial \xi}( G^i_k, \Xi^i_k ) - \sum_{j=1}^s a_{ji} \mu^j_k - b_i \mu_{k+1}, \delta \Xi^i_k \right\rangle \\
& \quad + \left\langle (d \tau^{-1}_{-h \Theta^i_k})^* \mu_k^i  + h b_i G^i_k \frac{\partial \ell}{\partial g}(G^i_k, \Xi^i_k ),  H^i_k \right\rangle  \\
& \quad + \left\langle  - (d \tau^{-1}_{h \xi_{k+1}})^* \mu_{k+1} + (d \tau^{-1}_{-h \xi_k})^* \mu_k  - \sum_{j=1}^s  (d \tau^{-1}_{h \Theta^j_k})^* \mu_k^j , \eta_k \right\rangle 
\end{align*}
Since $\mathbf{d} s_d(c_d) =0$ if and only if $\mathbf{d} s_d \cdot z = 0$ for all $z \in T_{c_d} \mathcal{C}_d$, one arrives at the desired equations with the elimination of $\mu_k^i$ and the introduction of the internal stage variables $M_k^i = \partial \ell/\partial \xi (G_k^i, \Xi_k^i) $ for $i=1,\cdots,s$.   Conversely, if $c_d$ satisfies \eqref{eq:lievprk1}--\eqref{eq:lievprk} then $\mathbf{d} s_d(c_d) =0$.

Consider the subset of $\mathcal{C}_d$ given by solutions of \eqref{eq:lievprk1}--\eqref{eq:lievprk}.  Let $\hat{s}_d$ denote the restriction of $s_d$ to this space. Since each of these solutions is determined by a point in $\mathcal{I}_{hp}$, one can identify this space with $\mathcal{I}_{hp}$, and hence, $\hat{s}_d : \mathcal{I}_{hp} \to \mathbb{R}$.   Since $\hat{s}_d$ is restricted to solution space,
\begin{align*}
\mathbf{d} &\hat{s}_d(g_0, \xi_0, \mu_0) \cdot (\delta g_0, \delta \xi_0, \delta \mu_0) = \left\langle (d \tau^{-1}_{- h \xi_N})^* \mu_N, g_N^{-1} \delta g_N \right\rangle - \left\langle  (d \tau^{-1}_{- h \xi_0})^* \mu_0,  g_0^{-1} \delta g_0 \right\rangle
  \text{.}
\end{align*}
Preservation of $\omega_{\mathcal{I}_{hp}}$ follows from $\mathbf{d}^2 \hat{s}_d = 0$.
\end{proof}
The external and internal stages of \eqref{eq:lievprk1}--\eqref{eq:lievprk} define update schemes on $G \times \mathfrak{g}^*$ and $(\mathfrak{g} \times \mathfrak{g} \times \mathfrak{g}^*)^s$, respectively.

\paragraph{St\"{o}rmer-Verlet Integrators on Lie Groups.}
The generalization of the St\"{o}rmer-Verlet method to Lie groups is given by evaluating \eqref{eq:lievprk1}--\eqref{eq:lievprk} at the following two-stage RK tableau (implicit trapezoidal rule),
\begin{center}
\begin{tabular}{c|cc}
0 &  $0$  & $0$ \\
 1  & $1/2$    & $1/2$ \\
\hline  
          & $1/2$    & $1/2$ 
\end{tabular} 
\end{center}
Given $h$ and $(g_k, \mu_k)$, the method determines $(g_{k+1}, \mu_{k+1})$ by solving the following system of equations:
\begin{align}
M^{1/2}_{k} &=   \frac{\partial \ell}{\partial \xi}(g_k, \Xi^{1}_k ), \label{eq:sv1} \\
M^{1/2}_k &= \frac{\partial \ell}{\partial \xi}(g_{k+1}, \Xi^{2}_k ), \label{eq:sv2} \\
 (d \tau^{-1}_{h \xi_{k+1}})^* M^{1/2}_{k} &=   (d \tau^{-1}_{-h \xi_k})^* \mu_k + \frac{h}{2}  g_k \frac{\partial \ell}{\partial g}(g_k, \Xi^1_k ), \label{eq:sv3} \\
 g_{k+1} & = g_k \tau \left( h \frac{1}{2} \left( \Xi_k^1 + \Xi_k^2 \right)  \right),\label{eq:sv4} \\
 \mu_{k+1} &=   M^{1/2}_k + \frac{h}{2} (d \tau_{-h \xi_{k+1}})^* g_{k+1} \frac{\partial \ell}{\partial g}(g_{k+1}, \Xi^2_k ) \text{.} \label{eq:sv5} 
\end{align}
In particular, one uses the following procedure:
\begin{itemize}
\item Eliminate $g_{k+1}$ in \eqref{eq:sv2} using \eqref{eq:sv4}.  Then solve for $M^{1/2}_k$, $\Xi^1_k$, and $\Xi^2_k$ using \eqref{eq:sv1}-\eqref{eq:sv3}.  This update is in general implicit.  
\item Update $g_{k+1}$ using \eqref{eq:sv4}.  This update is explicit.
\item Solve for $\mu_{k+1}$ using \eqref{eq:sv5}.  This update is explicit.
\end{itemize}
Observe that if the Lagrangian is separable, then \eqref{eq:sv3} is not implicit in the potential force term and one does not
need to eliminate $g_{k+1}$ in \eqref{eq:sv2} using \eqref{eq:sv4}.

\paragraph{Variational Euler on Lie Groups.} The variational Euler schemes come from evaluating \eqref{eq:lievprk1}--\eqref{eq:lievprk} with the following tableaus:
\begin{center}
\begin{tabular}{c|c}
0 &  0 \\
\hline
   & 1
\end{tabular} \text{,} \hspace{1.5in}
\begin{tabular}{c|c}
1 & 1 \\
\hline
   & 1
\end{tabular}  \text{.}\\
forward Euler \hspace{1.1in} backward Euler
\end{center}
The corresponding VPRK action sums take the following simple forms:
\[
\begin{array}{cc}
\begin{array}{c}
s_d^e = \sum_{k=0}^{N-1} h \left[ \ell(g_k, \xi_k) \right. \\
  + \left. \left\langle \mu_k, \tau^{-1}(g_k^{-1} g_{k+1})/h - \xi_k \right\rangle \right] \end{array} ~~~ &
  \begin{array}{c}
 s_d^i = \sum_{k=0}^{N-1} h \left[ \ell(g_{k+1}, \xi_{k+1}) \right. \\
  + \left. \left\langle \mu_{k+1}, \tau^{-1}(g_k^{-1} g_{k+1})/h  - \xi_{k+1} \right\rangle \right] \end{array}
 \end{array} 
\]
Given $h$ and $(g_k, \mu_k)$,  the forward variational Euler method determines $(g_{k+1}, \mu_{k+1})$ by solving the following system of equations:
\begin{equation}
\begin{cases}
\begin{array}{rll}
g_{k+1} &= & g_k \tau(h \xi_k) \text{,}\\
(d \tau^{-1}_{h \xi_{k+1}})^* \mu_{k+1} &= & (d \tau^{-1}_{-h \xi_k})^* \mu_k + h g_{k+1}       \frac{\partial \ell}{\partial g}(g_{k+1}, \xi_{k+1}) \text{,} \\ 
\mu_{k+1} &= &\frac{\partial \ell}{\partial \xi}(g_{k+1}, \xi_{k+1}) \text{.} 
\end{array}  
\end{cases}
\end{equation}
The backward variational Euler method determines $(g_{k+1}, \mu_{k+1})$ by solving the following system of equations:
\begin{equation}
\begin{cases}
\begin{array}{rllrll}
   g_{k+1} &= &g_k \tau(h \xi_{k+1})  \text{,} \\
(d \tau^{-1}_{h \xi_{k+1}})^* \mu_{k+1} &= & (d \tau^{-1}_{-h \xi_k})^* \mu_k + h g_k \frac{\partial \ell}{\partial g}(g_k, \xi_k) \text{,} \\ 
\mu_{k+1} &= &\frac{\partial \ell}{\partial \xi}(g_{k+1}, \xi_{k+1}) \text{.} 
\end{array}  
\end{cases}
\end{equation}

\paragraph{Euler-Poincar\'e Integrators.} In the case when the Lagrangian is $G$-left-invariant, the angular momentum updates in the above methods are identical and given by:
\begin{equation}
(d \tau^{-1}_{h \xi_{k+1}})^* \mu_{k+1} =  (d \tau^{-1}_{-h \xi_k})^* \mu_k  \label{eq:dlp}
\end{equation}

\paragraph{Examples.} We now give various examples of Euler-Poincar\'e integrators by making different choices of the map $\tau$ and evaluating \eqref{eq:dlp}.

\begin{description}

\item[(a) Matrix exponential.] Suppose 
\[ 
\tau = \exp(\xi),~~~\tau: \mathfrak{g} \to G \text{,} 
\]  
which is a local diffeomorphism.

Using standard convention the right trivialized tangent of the exponential map and its inverse are denoted by $\operatorname{dexp}: \mathfrak{g} \times \mathfrak{g} \to \mathfrak{g}$ and $\operatorname{dexp}^{-1}: \mathfrak{g} \times \mathfrak{g} \to \mathfrak{g}$, and are explicitly given by, 
\begin{equation}
\operatorname{dexp}(x) y = \sum_{j=0}^{\infty} \frac{1}{(j+1)!} \operatorname{ad}_x^j y,~~~
\operatorname{dexp}^{-1}(x) y = \sum_{j=0}^{\infty} \frac{B_j}{j!} \operatorname{ad}_x^j y \text{,} 
\label{eq:dexp}
\end{equation}   
where $B_j$ are the {\em Bernoulli numbers};   see \S3.4 of \cite{HaLuWa2006} for a detailed exposition and derivation.

Hence, \eqref{eq:dlp} takes the form, 
\begin{equation}
( \operatorname{dexp}^{-1}(h \xi_k))^* \mu_k =   (\operatorname{dexp}^{-1}(-h \xi_{k-1}))^* \mu_{k-1}    \label{eq:exphpvi} \text{.}
 \end{equation}

\item[(b) Pad\'e (1,1) approximant.] Suppose 
\begin{equation} 
\tau(\xi) = \operatorname{cay}(\xi) = (e -  \xi/2)^{-1} (e +\xi/2)  \label{eq:cay}  \text{,}
\end{equation} 
which is the Pad\'e (1,1) approximant to the matrix exponential and better known as the {\em Cayley transform}.   The Cayley transform maps to the group for quadratic Lie groups ($SO(n)$, the symplectic group $Sp(2 n)$, the Lorentz group $SO(3,1)$)  and the special Euclidean group $SE(3)$.

The right-trivialized tangent of the Cayley transform and its inverse are written below
\begin{equation} 
\operatorname{dcay}(x)y = (e - x/2)^{-1} y (e+x/2)^{-1},
~~~ \operatorname{dcay}^{-1}(x)y = (e-x/2) y (e + x/2) \label{eq:dcay} \text{.}  
\end{equation}  
For a derivation and exposition the reader is referred to \S 4.8.3 of \cite{HaLuWa2006}.   Using these expressions  \eqref{eq:dlp} can be written as,  
\begin{align}
\mu_k =&  \mu_{k-1} + \frac{h}{2}  \operatorname{ad}_{\xi_k}^* \mu_k  
  +  \frac{h}{2} \operatorname{ad}_{\xi_{k-1}}^*  \mu_{k-1}  \nonumber \\
&+ \frac{h^2}{4} \left( \xi_k^* \mu_k \xi_k^* - \xi_{k-1}^* \mu_{k-1} \xi_{k-1}^* \right) \label{eq:cayhpvi} \text{.}
\end{align}

 \item[(c) Pad\'e (1,0) or (0,1) approximant.] Rather than use the exact matrix exponential one can use a Pad\'e approximant, e.g., the Pad\'e (1,0) approximant 
\[
\exp(\xi) \approx e + \xi
\]  
or Pad\'e (0,1) approximant 
\[
\exp(\xi) \approx (e-\xi)^{-1} \text{.}
\] 
However, since a Pad\'e approximant is not guaranteed to lie on the group one needs to use a projector from $GL(n)$ to $G$.  In what follows $G=SO(n)$  will be considered where a natural choice of projector is given by skew symmetrization.

Suppose 
\[ 
\tau^{-1}(g) = \operatorname{skew}(g) = \frac{g-g^*}{2} \text{.} 
\] 
which comes from a first order approximant to the matrix exponential.  This map is a local diffeomorphism from a neighborhood of $e$ to a neighborhood of $0$ and its differential is the identity.   Its right trivialized tangent can be computed from its derivative: 
\begin{align*}
\operatorname{D}  \operatorname{skew}(g) \cdot \delta = \frac{\delta - \delta^*}{2} = \frac{ (\delta g^{-1} g ) - (\delta g^{-1} g)^* }{2} \text{.}
\end{align*} 
By definition of the right trivialized tangent of $\tau^{-1}$, it then follows that, 
\begin{equation} 
\operatorname{dskew}(x)(y)  = \frac{y \tau(x) - (y \tau(x))^*}{2}  \text{.}
\label{eq:dskew}
\end{equation} 
\cite{CaLe2003} obtained the following theorem that explicitly determines $\tau(\xi)$.  Moreover, they give necessary and sufficient conditions for its existence.

  \begin{theorem}
   Given $\xi \in \mathfrak{so(n)}$, a special orthogonal solution to the equation 
   \[
 \xi   = \frac{\tau(\xi)- \tau(\xi)^*}{2}
  \]  
  can be written as 
  \[ 
  \tau(\xi) = \xi + \left( \xi^2 + e \right)^{1/2} \text{,}
  \] 
  where  $ \left( \xi^2 + e \right)^{1/2}$ is a symmetric square root.  
  \end{theorem}

 \begin{proof} Since the skew-symmetric part of $g$ is $\xi$, one can write  $g$ as a sum of $\xi$ and a symmetric matrix $S$,  
 \[ 
 \tau(\xi) = S +  \xi \text{.} 
 \]
Observe that $\xi$ commutes with $\tau(\xi)$ since 
\[ 
2 \xi \tau(\xi) = ( \tau(\xi)-\tau(\xi)^* ) \tau(\xi)  = \tau(\xi)^2 - e = 2 \tau(\xi) \xi \text{.} 
\]
Moreover, $S$ satisfies an algebraic Riccati equation because, 
\[ 
\tau(\xi)^* \tau(\xi) = e \implies
S^2 +  S \xi - \xi S - (\xi^2 + e) = 0 \text{.}  
\] 
And since $\xi$ commutes with $S$ (because it commutes with $g$), 
\[ 
S^2 = ( \xi^2 + e) \text{,} 
\]
which completes the proof.  
 \end{proof}

  Hence, \eqref{eq:dlp} can be written as, 
  \begin{align}
& \frac{\mu_k   \left( h^2 \xi_k^2 + e \right)^{1/2} +  \left( h^2 \xi_k^2 + e \right)^{1/2} \mu_k }{2}  \nonumber \\
 & \qquad =   \frac{\mu_{k-1}   \left( h^2 \xi_{k-1}^2 + e \right)^{1/2} 
  +  \left( h^2 \xi_{k-1}^2 + e \right)^{1/2}  \mu_{k-1}}{2}  \nonumber \\
& \qquad \qquad + \frac{h}{2}  \operatorname{ad}_{\xi_k}^* \mu_k
    + \frac{h}{2}  \operatorname{ad}_{\xi_{k-1}}^* \mu_{k-1}   \label{eq:skewprojector}
  \end{align}
  
\end{description}

\section{Conclusion}

In this paper a left-trivialized Hamilton-Pontryagin principle is derived for mechanical systems on a Lie group $G$.   If the Lagrangian is left-invariant with respect to the action of $G$, it is shown that this left-trivialized HP principle unifies the Euler-Poincar\'e and Lie-Poisson descriptions.   In addition to its utility for implicit Lagrangian systems, the paper shows that this principle provides a practical way to design discrete Lagrangians.  In particular, the paper explains how one can discretize the kinematic constraint using a Runge-Kutta Munthe-Kaas (RKMK) method.  The paper shows that this leads to a novel generalization of variational partitioned Runge-Kutta methods from flat spaces to Lie groups.  In particular, one can generalize variational (or symplectic) Euler and St\"{o}rmer-Verlet methods to Lie groups in this fashion.  These methods inherit many of their attractive properties on flat spaces: efficiency, order of accuracy, symplecticity, symmetry, etc.

Part II of this paper will develop a basic numerical analysis of these methods and report on numerical experiments on a class of nonreversible mechanical systems on Lie groups (moving rigid body systems) and chaotic dynamics of an underwater vehicle.   To be specific the paper will:
\begin{itemize}
\item prove order of accuracy of the VPRK integrators presented in this paper by invoking the variational proof of order of accuracy \cite{MaWe2001};
\item explain the numerics behind the Poincar\'e sections provided in Figure~\ref{fig:uv};
\item demonstrate the superiority of these VPRK integrators compared to symmetric rigid body integrators  when applied to a nonreversible system such as a rigid body on a turntable.  
\end{itemize}

\newpage



\begin{thebibliography}{}


\bibitem[Austin et al.(1993)]{AuKrWa1993}
Austin, M.~A., P.~S.~Krishnaprasad, and L.~Wang [1993],
Almost-Poisson Integration of Rigid Body Systems.
{\em Journal of Computational Physics}, {\bf 107},  105--117.

\bibitem[Bobenko and Suris(1999a)]{BoSu1999a}
Bobenko, A.~I., and Y.~B.~Suris [1999a], Discrete Lagrangian reduction,  
discrete Euler-Poincar\'e equations, and semi-direct products. 
 {\em Lett.~Math.~Phys.}, {\bf 49}, 79--93.

\bibitem[Bobenko \& Suris(1999b)]{BoSu1999b}
Bobenko, A.~I., and Y.~B.~Suris [1999b],  Discrete time Lagrangian Mechanics on Lie groups, 
with an application to the Lagrange top. {\em Comm.~Math.~Phys.}, {\bf 204},  147--188.

\bibitem[Bou-Rabee(2007)]{Bou-Rabee2007}
Bou-Rabee, N. [2007], {\em Hamilton-Pontryagin integrators on Lie groups}, PhD
  thesis, California Institute of Technology.

\bibitem[Bou-Rabee \& Owhadi(2007b)]{BoOw2007b}
Bou-Rabee, N., and H.~Owhadi [2007b],  Stochastic Variational Partitioned Runge-Kutta 
Integrators for Constrained Systems. Submitted; arXiv:0709.2222.

\bibitem[Cardoso \& Leite(2003)]{CaLe2003}
Cardoso, J., R., and F.~Leite [2003],  
The Moser-Veselov Equation.  
{\em Linear Algebra and its Applications}, {\bf 360}, 237--248.

\bibitem[Celledoni  \& Iserles(2001)]{CeIs2001}
Celledoni, E., and A.~Iserles [2001],
Methods for the approximation of the matrix exponential in a Lie-algebraic setting.
{\em IMA J.~Num.~Anal.}, {\bf 21},  463--488.

\bibitem[Cendra et~al.(2003)Cendra, Marsden, Pekarsky, and Ratiu]{CeMaPeRa2003} Cendra, H., J.~E. Marsden, S.~Pekarsky, and T.~S. Ratiu [2003], Variational principles for {L}ie-{P}oisson and {H}amilton-{P}oincar{\'e} equations, {\em Mosc. Math. J.}, \textbf{3}, 833--867.

\bibitem[Dahlquist(1975)]{Da1975}
Dahlquist, G.~[1975],
Error analysis for a class of methods for stiff nonlinear initial boundary value problems.
{\em Lecture Notes in Mathematics}, {\bf 506}, 60--74.

\bibitem[Dullweber(1997)]{DuLeMc1997}
Dullweber, A., B.~Leimkuhler, and R.~McLachlan [1997],
Symplectic splitting methods for rigid body molecular dynamics.
{\em J.~Chem.~Phys.}, {\bf 107},  5840--5851.


\bibitem[Eng\o~\& Faltinsen(2001)]{EnFa2001}
Eng\o, K., and  S.~Faltinsen [2001],
Numerical integration of Lie-Poisson systems while preserving 
coadjoint orbits and energy.  {\em SIAM J.~Numer.~Anal.}, 
{\bf 39}, 128--145.


\bibitem[Feng(1986)]{Fe1986}
Feng, K.~[1986], 
Difference Schemes for Hamiltonian Formalism and Symplectic Geometry.
{\em J.~Comp.~Math.}, {\bf 4}, 279--289.


\bibitem[Gonzalez \& Simo(1996)]{GoSi1996}
Gonzalez, O., and J.~C.~Simo [1996],
On the stability of symplectic and energy-momentum algorithms for nonlinear
Hamiltonian systems with symmetry.
{\em Comput.~Methods Appl.~Mech.~Eng.}, {\bf 134},  197--222.


\bibitem[Hairer et~al.(2006)Hairer, Lubich, and Wanner]{HaLuWa2006}
Hairer, E., C.~Lubich, and G.~Wanner [2006], {\em Geometric numerical
  integration}, volume~31 of {\em Springer Series in Computational
  Mathematics}. Springer-Verlag, Berlin, second edition.

\bibitem[Holm et~al.(1998)]{HoMaRa1998}
Holm, D.~D., J.~E. Marsden, and T.~S. Ratiu [1998], The
  {E}uler--{P}oincar{\'{e}} equations and semidirect products with applications
  to continuum theories, {\em Adv. in Math.} {\bf 137}, 1--81.
  
\bibitem[Holmes et al.(1998)]{HoJeLe1998}
Holmes, P., J.~Jenkins, and N.~E.~Leonard [1998],
Dynamics of the Kirchhoff equations I: Coincident centers of
gravity and buoyancy.  {\em Phys.~D}, {\bf 118}, 311--342.

\bibitem[Iserles et~al.(2000)Iserles, Munthe-Kaas, N{\o}rsett, and
  Zanna]{IsMuNoZa2000} Iserles, A., H.~Z. Munthe-Kaas, S.~P. N{\o}rsett, and A.~Zanna [2000], Lie-group methods. {\em Acta numerica,}, {\bf 9}, 215--365. 

\bibitem[Iserles et al.(2001)]{Is2001}
Iserles, A. [2001],
On Cayley-Transform methods for the discretization of
Lie-group equations.  {\em Found.~Comp.~Maths}, {\bf 1}, 129--160.


\bibitem[Kane et al.(2000)]{KaMaOrWe2000}
Kane, C., J.~E.~Marsden, M.~Ortiz, and M.~West [2000],
Variational integrators and the Newmark algorithm for conservative
and dissipative mechanical systems.
{\em Int.~J.~Num.~Meth.~Eng'g.}, {\bf 49}, 1295--1325.

\bibitem[Kane et al.(1999)]{KaMaOr1999}
Kane, C., J.~E.~Marsden, and M.~Ortiz [1999],
Symplectic-Energy-Momentum Preserving Variational Integrators.
{\em J.~Math.~Phys.}, {\bf 40}, 3353--3371.

\bibitem[Kanso et al.(2005)]{KaMaRoMe2005}
Kanso, E., J.~E.~Marsden, C.~W.~Rowley, and J.~B.~Melli-Huber [2005],
Locomotion of articulated bodies in a perfect fluid. 
{\em J.~Nonlinear Sci.}, {\bf 15}, 255--289.

\bibitem[Kharevych et al.(2006)]{KWTKMSD2006}
Kharevych, L., Weiwei, Y.~Tong, E.~Kanso, J.~E.~Marsden, P.~Schroder, and
M.~Desbrun [2006],  Geometric, Variational Integrators for
Computer Animation. {\em Eurographics/ACM SIGGRAPH Symposium
on Computer Animation.  }


\bibitem[Lall \& West(2006)]{LaWe2006}
Lall, S.~and M.~West [2006],
{\em Discrete variational Hamiltonian mechanics.} 
Journal of Physics A: Mathematical and General, \textbf{39}, 5509--5519.

\bibitem[Leimkuhler \& Reich(2004)]{LeRe2004}
Leimkuhler, B., and S.~Reich [2004], 
Simulating Hamiltonian Dynamics. {\em Cambridge Monographs on 
Applied and Computational Mathematics}, {\bf 14}. 

\bibitem[Leok et al.(2005)Leok, McClamroch , and Lee]{LeMcLe2005} 
Leok, M., N.~H.~McClamroch, and  T.~Lee [2005],
A Lie Group Variational Integrator for the Attitude Dynamics 
of a Rigid Body with Applications to the 3D Pendulum.
{\em Proc.~IEEE Conf.~on Control Applications}, 962--967.

\bibitem[Lee et~al.(2007)Lee, Leok, and McClamroch]{LeLeMc2007}
Lee, T., M.~Leok, and N.~H. McClamroch [2007], Lie group variational
  integrators for the full body problem, {\em Comput. Methods Appl. Mech.
  Engrg.} \textbf{196}, 2907--2924.

\bibitem[Lew et al.(2003)]{LeMaOrWe2003}  
Lew, A., J.~E.~Marsden, M.~Ortiz, and M.~West [2003],
Asynchronous Variational Integrators.  
{\em Arch.~Rational Mech.~Anal.}, {\bf 167}, 85--146.

\bibitem[Lew et al.(2004a)]{LeMaOrWe2004a}
Lew, A., J.~E.~Marsden, M.~Ortiz, and M.~West [2004a],
An Overview of Variational Integrators.
{\em Finite Element Methods: 1970's and Beyond}, CIMNE Barcelona, 98--115.

\bibitem[Lew et~al.(2004)Lew, Marsden, Ortiz, and West]{LeMaOrWe2004b} Lew, A., J.~E. Marsden, M.~Ortiz, and M.~West [2004], Variational time integrators, {\em Internat. J. Numer. Methods Engrg.} \textbf{60}, 153--212.

\bibitem[Lewis \& Simo(1994)]{LeSi1994}
Lewis, D.~and J.~C.~Simo [1994],  
Conserving Algorithms for the 
Dynamics of Hamiltonian Systems on Lie Groups.
{\em J.~Nonlinear Science}, {\bf 4}, 253--299.

\bibitem[Lewis \& Simo(1996)]{LeSi1996}
Lewis, D.~and J.~C.~Simo [1996],  
Conserving Algorithms for the N-dimensional rigid body.
{\em Fields.~Inst.~Comm.}, {\bf 10}, 121--139.

\bibitem[Leyendecker et~al.(2007)Leyendecker, Ober-Bl\"{o}baum, Marsden, and
  Ortiz]{LeObMaOr2007}
Leyendecker, S., S.~Ober-Bl\"{o}baum, J.~E.~Marsden, and M.~Ortiz [2007], Discrete
  mechanics and optimal control for constrained multibody dynamics.
\newblock In {\em {Proceedings of the 6th International Conference on Multibody
  Systems, Nonlinear Dynamics, and Control, ASME International Design
  Engineering Technical Conferences},}, volume Proceedings of IDETC/MSNDC 2007,
  pages 1--10.


\bibitem[Leyendecker, Marsden, and Ortiz(2007)]{LeMaOr2007}
Leyendecker, S., J.~E.~Marsden, and M.~Ortiz [2007], Variational integrators for
  constrained mechanical systems, {\em (preprint)}.


\bibitem[Livens(1919)]{Livens1919}
Livens, G.~H. [1919], On Hamilton's principle and the modified function in
  analytical dynamics, {\em Proc. Roy. Soc. Edingburgh} \textbf{39}, 113.



\bibitem[MacIver et al.(2004)]{MaFoBu2004}
MacIver, A., E.~Fontaine, and J.~W.~Burdick [2004],
Designing future underwater vehicles: 
Principles and mechanisms of the weakly electric fish.
{\em IEEE J. Oceanic Engineering}, {\bf 29}, 651--659.

\bibitem[Mackay(1992)]{Ma1992}
MacKay, R.~[1992],  Some Aspects of the Dynamics of
Hamiltonian Systems. {\em The Dynamics of Numerics
and the Numerics of Dynamics},  Clarendon Press, 137--193.

\bibitem[Marsden \& Scheurle(1993)]{MaSc1993}
Marsden, J.~E., and J.~Scheurle [1993],  
The reduced Euler-Lagrange Equations.
{\em Fields.~Inst.~Comm.}, {\bf 1}, 139--164.

\bibitem[Marsden, Pekarsky, and Shkoller(1998)]{MaPeSh1998}
Marsden, J.~E., S.~Pekarsky, and S.~Shkoller [1998], 
Discrete Euler-Poincar\'e and Lie-Poisson Equations.
{\em Nonlinearity}, {\bf 12}, 1647--1662.

\bibitem[Marsden \& Ratiu(1999)]{MaRa1999}
Marsden, J.~E., and T.~Ratiu [1999],
 {\it Introduction to Mechanics and Symmetry}. Springer 
Texts in Applied Mathematics, {\bf 17}. 

\bibitem[Marsden \& West(2001)]{MaWe2001}
Marsden, J.~E., and M.~West [2001],  
Discrete mechanics and variational integrators. {\em Acta Numerica},
{\bf 10}, 357--514.

\bibitem[McLachlan \& Scovel(1995)]{McSc1995}
Mclachlan, R., and C.~Scovel [1995],  
Equivariant constrained symplectic integration. 
{\em J.~Nonlinear Science},
{\bf 5}, 233--256.

\bibitem[Moser \& Veselov(1991)]{MoVe1991}
Moser, J., and A.~P.~Veselov [1991],  
Discrete versions of some classical integrable systems
and factorization of matrix polynomials. 
{\em Comm.~Math.~Phys.},
{\bf 139}, 217--243.

\bibitem[Munthe-Kaas(1995)]{Mu1995}
Munthe-Kaas, H. [1995],
Lie-Butcher theory for Runge-Kutta methods.
{\em BIT}, 572--587.

\bibitem[Munthe-Kaas(1998)]{Mu1998}
Munthe-Kaas, H. [1998],
\newblock Runge-Kutta methods on lie groups.
{\em BIT}, {\bf 38}, 92--111.

\bibitem[Munthe-Kaas \& Owren(1999)]{MuOw1999}
Munthe-Kaas, H. and B.~Owren [1999],
Computations in a free lie algebra.
{\em Phil. Trans Royal Society A}, {\bf 357}, 957--982.

\bibitem[Munthe-Kaas \& Zanna(1997)]{MuZa1997}
Munthe-Kaas, H. and A.~Zanna [1997],
Numerical integration of differential equations on homogeneous manifolds. 
{\em Springer}, 305--315.



\bibitem[Ruth(1983)]{Ru1983}
Ruth, R.~D.~[1983], 
A canonical integration technique.
{\em IEEE Transactions on Nuclear Science}, {\bf 30}, 2669--2671.


\bibitem[Suris(1990)]{Su1990}
Suris, Y.~B.~[1990],  
Hamiltonian methods of Runge-Kutta type and their variational 
interpretation.  {\em Math.~Modelling},  {\bf 2},  78--87.

\bibitem[Simo \& Vu-Quoc(1988)]{SiVu1988}
Simo, J.~C.~and L.~Vu-Quoc [1988],
On the Dynamics in Space of Rods Undergoing Large 
Motions - A Geometrically Exact Approach.
{\em Comp.~Methods in Appled Mech.~and Eng'g.}, {\bf 66}, 125--161.

\bibitem[Simo \& Wong(1991)]{SiWo1991}
Simo, J.~C.~and T.~S.~Wong [1991],  
Unconditionally Stable Algorithms for Rigid Body 
Dynamics That Exactly Preserve Energy and Momentum.
{\em Int.~J.~Numer.~Methods Engrg.}, {\bf 31}, 19--52.


\bibitem[Veselov(1988)]{Ve1988}
Veselov, M.~[1988],  Integrable discrete-time systems and difference operators.
{\em Func.~An.~and Appl.}, {\bf 22}, 83--94.

\bibitem[de Vogelaere(1956)]{Vo1956}
de Vogelaere, R.~[1956],
Methods of Integration which Preserve 
the Contact Transformation Property of the Hamiltonian Equations.
Report No.~4, Dept.~Math., Univ.~of Notre Dame.


\bibitem[Wendlandt and Marsden(1997)]{WeMa1997}
Wendlandt, J.~M. and J.~E. Marsden [1997], Mechanical integrators derived from
  a discrete variational principle, {\em Physica D} \textbf{106}, 223--246.


\bibitem[Yoshimura and Marsden(2006a)]{YoMa2006a}
Yoshimura, H. and J.~E. Marsden [2006], Dirac structures and Lagrangian
  Mechanics Part I: Implicit Lagrangian systems, {\em J. Geom. and Physics}
  \textbf{57}, 133--156.

\bibitem[Yoshimura and Marsden(2006b)]{YoMa2006b}
Yoshimura, H. and J.~E. Marsden [2006b], Dirac structures in Lagrangian
  mechanics Part II: Variational structures, {\em J. Geom. and Physics}
  \textbf{57}, 209--250.

\end{thebibliography}
\end{document}